\documentclass[12pt]{amsart}
\usepackage{amsmath,amssymb,amscd,array, mathrsfs }

\usepackage{amsmath,amscd,amssymb,amsthm,array}
\usepackage{color}

\usepackage{amsmath,amssymb,mathrsfs,amsthm, tikz-cd,mathrsfs}
\usepackage{comment}
\def\url#1{\expandafter\s

\tring\csname #1\endcsname}

\def\mmat #1,#2,#3,#4,{\text{\small\arraycolsep=3pt $
\begin{pmatrix}#1&#2\\#3&#4\end{pmatrix}$}}

\usepackage{hyperref}
\usepackage{capt-of}

\usepackage{multirow}
\usepackage[all]{xy}

\usepackage{comments}
\newComments\SBe{Said}{blue}
\newComments\SBo{Sofiane}{blue}
\newComments\AM{Nacer}{blue}
\newComments\DL{DL}{red}
\newComments\QEh{QEh}{blue}

\def\mmat #1,#2,#3,#4,{\text{\small\arraycolsep=3pt $
\begin{pmatrix}#1&#2\\#3&#4\end{pmatrix}$}}

\usepackage{lscape}
\usepackage{tikz-cd}
\usepackage{enumerate}

\usepackage{DLdef1}
\usepackage{multicol}

\def\mmat #1,#2,#3,#4,{\text{\small\arraycolsep=3pt $
\begin{pmatrix}#1&#2\\#3&#4\end{pmatrix}$}}

\renewcommand {\ssbegin}[2][*]
 {\refstepcounter{subsection}%
\if#1*
\addcontentsline{toc}{subsection}{\thesubsection.\hspace*{1pc} #2}%
\else
\addcontentsline{toc}{subsection}{\thesubsection.\hskip 1pc #2. #1}%
\fi
 \def \secno {\gdef \secno {}{\ssecfont
\thesubsection.\hskip 2ex}%
 }%
 \begin{#2}}

\renewcommand {\sssbegin}[2][*]
 {\refstepcounter{subsubsection}
\if#1*
\addcontentsline{toc}{subsubsection}{\thesubsubsection.\hskip 1pc #2}%
\else
\addcontentsline{toc}{subsubsection}{\thesubsubsection.\hskip 1pc #2. #1}
\fi
 \def \secno {\gdef \secno {}{\ssecfont \thesubsubsection.\hskip 2ex}%
 }%
 \begin{#2}}

\renewcommand {\parbegin}[2][*]
 {\refstepcounter{paragraph}
\if#1*
\addcontentsline{toc}{paragraph}{\theparagraph.\hskip 1pc #2}%
\else
\addcontentsline{toc}{paragraph}{\theparagraph.\hskip 1pc #2. #1}
\fi
 \def \secno {\gdef \secno {}{\ssecfont \theparagraph.\hskip 2ex}%
 }%
 \begin{#2}}

\rmnameii{etr}{etr}
\rmnameii{evv}{ev}

\DeclareMathOperator{\K}{\mathbb{K}}
\newcommand {\A}{{\cal{A}}}

\newcommand{\B}{\rm{B}}

\newcommand{\Ll}{{\mathrm{L}}}
\newcommand{\Rr}{{\mathrm{R}}}

\newcommand{\al}{\alpha}

\newcommand{\la}{\lambda}

\newcommand{\In}{ \{\cdot,\cdot\}}

\def\br{[\cdot,\cdot]}
\newcommand{\esp}{\quad\mbox{and}\quad}
\begin{document}

\title[Anti-quadratic nearly associative algebras ]{Anti-quadratic nearly associative algebras}

\author{Sofiane Bouarroudj}
\thanks{{\bf Corresponding author}: Sofiane Bouarroudj.\\
S.~B. was supported by the grant NYUAD-065.}

\address{Division of Science and Mathematics, New York University Abu Dhabi, P.O. Box 129188, Abu Dhabi, United Arab Emirates.}
\email{sofiane.bouarroudj@nyu.edu}

\author{Hamza El Ouali}
\address{Division of Science and Mathematics, New York University Abu Dhabi, P.O. Box 129188, Abu Dhabi, United Arab Emirates.}
\email{eloualihamza11@gmail.com}


\keywords{Nearly associative algebras,  anti-$\Ll$-algebras, anti-$\Rr$-algebras, anti-quadratic forms,   $TT^*$-extension, double extension}

 \subjclass[2020]{17A30; 17A60; 17D25}

\begin{abstract}
An \emph{anti-quadratic} nearly associative algebra $(A, \bullet, \omega)$ is a nearly associative algebra equipped with a non-degenerate skew-symmetric bilinear form $\omega$ satisfying the anti-associativity condition $\omega(u \bullet v, w) = -\omega(u, v \bullet w)$ for all  $u, v, w \in A.$   In this paper, we investigate the structural properties and extension theory  of these algebras. We prove that any anti-quadratic nearly associative algebra is inherently nilpotent of index at most $6$. Furthermore, we introduce a novel $T T^*$-extension method and show that every anti-quadratic nearly associative algebra arises from this construction. We also introduce the notion of a double extension and prove that every such algebra can be constructed by iterated double extensions. Finally, as an application, we give a complete classification of anti-quadratic nearly associative algebras of dimension at most $4$.

\end{abstract}


\maketitle

\thispagestyle{empty}
\setcounter{tocdepth}{2}
\section{Introduction}

An algebra $(\mathcal{A}, \bullet)$ over a field $\mathbb{K}$ is said to be \textit{nearly associative} if it satisfies the following trilinear identity:
\[
(u \bullet v) \bullet w = v \bullet (w \bullet u), \quad \text{for all } u,v,w \in \mathcal{A}.
\]
Originally introduced by Dassoundo and Silvestrov in \cite{DS}, this family of non-associative structures is distinguished by a ternary relation that permutes the elements cyclically. Consequently, nearly associative algebras can be framed as a specific instance of the broader algebraic classes investigated by Kosier and Osborn in \cite{Osborn}. A systematic study of their structural and geometric properties was subsequently carried out in \cite{BBR}. Among other results, it was shown that the Wedderburn Principal Theorem holds for this class of algebras. These structures were also investigated under the name of \emph{shift associative algebras} in \cite{Ivan}, where the authors obtained a complete classification in dimensions at most four. Further insights on their algebraic features can be found in \cite{BBR1}, which develops the bialgebra framework for both nearly associative and LR-algebras. 

Given any arbitrary non-associative algebra $(\mathcal{A}, \bullet)$, one can naturally associate two distinct structures by defining the anti-commutative and the commutative products, respectively given by:
\[
[u,v] := u \bullet v - v \bullet u,
\qquad
\{u,v\} := u \bullet v + v \bullet u, \quad \text{for all } u,v \in \mathcal{A}.
\]
These operations endow the underlying space with the structures of an anti-commutative algebra $\mathcal{A}^- = (\mathcal{A}, \br)$ and a commutative algebra $\mathcal{A}^+ = (\mathcal{A}, \In)$. Recall that $(\mathcal{A}, \bullet)$ is termed \textit{Lie-admissible} (resp. \textit{Jordan-admissible} or \textit{Jacobi-Jordan-admissible}) if $\mathcal{A}^-$ constitutes a Lie algebra (resp. if $\mathcal{A}^+$ forms a Jordan or Jacobi-Jordan algebra). A remarkable feature established in \cite{BBR} is that every nearly associative algebra is simultaneously Lie-admissible and Jordan-admissible, thereby embedding this class within the intersection of these two classical algebraic theories.

In this paper, we introduce and investigate the structural properties of \emph{anti-quadratic nearly associative algebras}. Let us first clarify the foundational definition of this geometric structure. An anti-quadratic algebra is a triple $(A, \bullet, \omega)$ consisting of an algebra $(A, \bullet)$ equipped with a non-degenerate, skew-symmetric bilinear form $\omega$ satisfying the anti-associativity condition:
\begin{equation}\label{anti-assoc}
\omega(u \bullet v, w) = -\omega(u, v \bullet w), \quad \text{for all } u, v, w \in A.
\end{equation}

It is worth noting that in certain contexts, particularly within the study of $2$-Engel groups, these structures have traditionally appeared under the terminology of \emph{symplectic alternating algebras} (see, e.g., \cite{T}). Some general theory was developed  in \cite{TTT, T1}, and their structure theory in the nilpotent case was extensively explored by Sorkatti and Traustason in \cite{Lyla1}. Furthermore, a complete classification of such $10$-dimensional nilpotent structures over any field was established in \cite{Layla3}. 

More recently, in \cite{BE1}, anti-quadratic Lie algebras were investigated under the terminology of $T^*$-presymplectic Lie algebras. In that work, the authors introduced the concept of $TT^*$-extensions for $2$-step nilpotent Lie algebras equipped with a non-degenerate, skew-symmetric bilinear form. This construction serves as a natural analogue to the classical $T^*$-extension method for non-associative algebras endowed with a non-degenerate, \emph{symmetric} bilinear form, a seminal framework originally introduced by Bordemann in \cite{Bordemann} for the study of metrised algebraic structures.

In our current framework, however, we deliberately adopt the terminology of \emph{anti-quadratic algebras} rather than \emph{symplectic alternating algebras}. This choice is motivated by a deep and natural algebraic analogy. Recall that a standard \emph{quadratic algebra} is defined by a non-degenerate, symmetric bilinear form that is associative, meaning $B(u \bullet v, w) = B(u, v \bullet w)$ for all $u,v,w \in A$. Replacing the symmetry with skew-symmetry and the associativity with the anti-associativity identity given in \eqref{anti-assoc} we arrive at our framework. Thus, the term \emph{anti-quadratic} serves as a more coherent structural counterpart to quadratic geometry, while concurrently preventing any potential confusion with other widely established symplectic structures in classical differential geometry and Lie theory.

\subsection{Outline of the paper}
Section \ref{back} collects the preliminary material required for our work. For completeness and to establish notation, we include a review of key concepts, some of which may be familiar to specialists.

In Section \ref{S3}, we investigate the structural and representation properties of nearly associative algebras. We begin by showing that, within the anti-quadratic framework, the notions of nearly associative algebras, anti-$\Ll$-algebras, and anti-$\Rr$-algebras coincide completely, see Theorem \ref{cara-nearly}. Following this characterization, we prove that any anti-quadratic nearly associative algebra is inherently nilpotent of index at most $6$, see Theorem \ref{nil}. Subsequently, we turn our attention to their representation theory by introducing both adjoint and coadjoint representations, and establishing the exact conditions required for the existence of the latter. Finally, we provide a duality result, demonstrating that a nearly associative algebra admits an anti-quadratic bilinear form if and only if its adjoint and coadjoint representations are isomorphic as $A$-bimodules, See Theorem \ref{antieqdual}.

In Section \ref{S4}, we introduce a novel construction inspired by the classical $T T^*$-extension method, which was originally developed by Bordemann \cite{Bordemann} for quadratic non-associative algebras. This method was later adapted to anti-quadratic Lie algebras in \cite{BE1}, and is referred to as a Lagrangian extension for symplectic Lie algebras in \cite{BC} (see also \cite{BM, BR} for the corresponding notion in the super setting). Our goal is to extend the core philosophy of this construction to the setting of anti-quadratic nearly associative algebras, taking into account the structural peculiarities of the dual maps and the underlying anti-quadratic bilinear forms. Remarkably, we prove that this construction is exhaustive; that is, we show that every anti-quadratic nearly associative algebra can be obtained via a $T T^*$-extension, see Theorem \ref{Converse-T}.

In Section \ref{S5}, we introduce the notions of annihilator extensions and almost semidirect products for nearly associative algebras. These two structural tools play a crucial r\^ole in our framework, as they constitute the foundational building blocks required to define and construct the double extension of anti-quadratic nearly associative algebras, see Theorem \ref{construction-ex}. This construction, inspired by the classical double extension method of Medina and Revoy \cite{Medina-Revoy} for quadratic Lie algebras, provides an effective inductive way to generate new anti-quadratic structures from lower-dimensional ones. Finally, we prove that every anti-quadratic nearly associative algebra can be obtained by means of such a double extension, see Theorem \ref{Th-ex1} and Corollary \ref{succession-double}.

In Section~\ref{S6}, as an application of the double extension theory, we obtain a complete classification of anti-quadratic nearly associative algebras of dimension at most 4.

\section{Backgrounds}\label{back}
Throughout this work, all vector spaces are assumed to be finite-dimensional over a field $\mathbb{K}$ of characteristic zero.

We begin by recalling several fundamental definitions and basic concepts that will be used throughout the paper.
\subsection{Symplectic vector spaces}
Let $E$ be a vector space. For elements \(x,y\in E\), we denote their duals by \(x^*,y^*\in E^*\). For every \(x^*\otimes y^*\in E^*\otimes E^*\) and  \(u\otimes v\in E\otimes E\), we define the pairing by 
\[
\langle x^*\otimes y^*, u \otimes v\rangle :=x^*(u)\,y^*(v).
\]
The wedge product of $x^*$ and $y^*$ is given by
\[
x^*\wedge y^*:=x^*\otimes y^* -y^*\otimes x^*.
\]

A \emph{symplectic} vector space $(E, \omega)$ is a  vector space $E$ endowed with a nondegenerate, skew-symmetric bilinear form $\omega$. A \emph{symplectic basis} of $(E, \omega)$ is a basis $\{e_1, \dots, e_{2n}\}$ of $E$ such that 
\[
\omega(e_{2i-1}, e_{2i}) = -\omega(e_{2i}, e_{2i-1}) = 1, \qquad \text{for } i = 1, \dots, n,
\]
with all other pairings between basis elements being zero.

Let $F$ be a subspace of $(E, \omega)$. The orthogonal complement of $F$ with respect to $\omega$ is the subspace defined by 
$$
F^{\perp} = \{x \in E \mid \omega(x, y) = 0 \text{ for all } y \in F\}.
$$ The subspace $F$ is said to be \emph{degenerate} (resp., \emph{nondegenerate}, \emph{totally isotropic}, \emph{Lagrangian}) if $F \cap F^{\perp} \neq \{0\}$ (resp., $F \cap F^{\perp} = \{0\}$, $F \subset F^{\perp}$, $F = F^{\perp}$). 

Given an endomorphism $f : E \rightarrow E$, its adjoint $f^*$ with respect to $\omega$ is the endomorphism of $E$ uniquely defined by 
$$
\omega(f(x), y) = \omega(x, f^*(y)) \text{ for all $x, y \in E$}.
$$
\subsection{Jordan and Jacobi-Jordan algebras}
An algebra $(\mathcal{J}, \{\cdot,\cdot\})$ is called a \emph{Jordan algebra} if, for all $u,v \in \mathcal{J}$, the following identities hold:
\begin{align*}
\{u,v\} &= \{v,u\}, \\
\{\{u,v\},u^2\} &= \{u,\{v,u^2\}\}, \qquad \text{(Jordan identity)}
\end{align*}
where $u^2 := \{u,u\}$.

An algebra $(\mathcal{J}, \In)$ is called a \emph{Jacobi-Jordan algebra} if, for all $u,v,w \in \mathcal{J}$, the following identities  hold:
\begin{align*}
\{u,v\} &= \{v,u\}, \\
\{\{u,v\},w\} + \{\{v,w\},u\} + \{\{w,u\},v\} &= 0, \qquad \text{(Jacobi identity)}.
\end{align*}
This class of non-associative algebras constitutes a subclass of commutative Jordan algebras. Furthermore, Jacobi-Jordan algebras are well known to be nilpotent (see, e.g., \cite{BF, OK}).

Let $(A,\bullet )$ be a non-associative algebra. On the underlying vector space of $A$, we can define the commutator $\br$ and the anti-commutator $\In$ respectively by
\begin{equation}\label{skew-sym}
[u,v] := u\bullet v - v\bullet u, 
\qquad 
\{u,v\} := u\bullet v + v\bullet u,
\end{equation}
for all $u,v \in A$. 

We denote by $A^{-} := (A, \br)$ the algebra endowed with the commutator product, and by $A^{+} := (A, \In)$ the one endowed with the anti-commutator product.

\sssbegin{Definition}
Let $(A, \bullet)$ be a non-associative algebra. Then:
\begin{enumerate}
    \item $(A, \bullet)$ is said to be \emph{Lie-admissible} if the algebra $A^{-}$ is a Lie algebra;
    \item $(A, \bullet)$ is said to be \emph{Jordan-admissible} if the algebra $A^{+}$ is a Jordan algebra;
    \item $(A, \bullet)$ is said to be \emph{Jacobi–Jordan-admissible} if the algebra $A^{+}$ is a Jacobi-Jordan algebra.
\end{enumerate}
\end{Definition}
Let $(A, \bullet)$ be a non-associative algebra. As usual, we denote by $\Ll, \Rr : A \to \mathrm{End}(A)$ the left and right multiplication operators, respectively. More precisely, for each $u \in A$, these operators are defined by
\[
\Ll_u(v) := u \bullet v, \qquad \Rr_u(v) := v \bullet u,
\]
for any $v \in A$. The two-sided ideal $A^2 = A \bullet A$ is defined as the span of all products:
\[
A \bullet A := \mathrm{span}\{u \bullet v \mid u, v \in A\},
\]
and the annihilator of $A$, denoted by $N(A)$, is given by
\[
N(A) := \{u \in A \mid \Ll_u = \Rr_u = 0\}.
\]

Recall that the \emph{derived series} of an algebra $(A, \bullet)$ is defined by
\[
A \supseteq A^{(1)} \supseteq A^{(2)} \supseteq \cdots,
\]
where $A^{(1)} := A^2$ and $A^{(i+1)} := (A^{(i)})^2$ for all $i \geq 1$.  
The algebra $A$ is said to be \emph{solvable} if there exists an integer $r \geq 1$ such that $A^{(r)} = \{0\}$.

An algebra $(A, \bullet)$ is called \emph{nilpotent} if there exists an integer $r$ such that for all $u_1,\ldots, u_r \in A$, the product $u_1\bullet \cdots \bullet u_r=0$, regardless of the parenthesization.
 \subsection{Pre-symplectic algebras}
Let $(A, \bullet)$ be a non-associative algebra. A bilinear form $\omega : A \times A \longrightarrow \mathbb{K}$ is said to be:
\begin{enumerate}
    \item \textit{skew-symmetric} if 
    $
    \omega(u, v) = -\omega(v, u),$ for all  $u, v \in A$;
    
    \item \textit{nondegenerate} if \(\omega(u,v)=0\) for all \(v\in A\) implies that \(u=0\);
    
    \item \emph{anti-invariant \textup{(}or anti-associative \textup{)}} if
    $
    \omega(u \bullet v, w) = -\omega(u, v \bullet w),$ for all $u, v, w \in A;
    $

\end{enumerate}

A  non-associative algebra $(A, \bullet)$  endowed with a nondegenerate skew-symmetric bilinear form $\omega$. The triple $(A, \bullet, \omega)$ is said to be a \emph{pre-symplectic non-associative algebra} and the bilinear form is referred to as a \emph{pre-symplectic structure} on $A$.

\sssbegin{Definition}[anti-quadartic non-associative algebra]\label{left-invariant pseudo-Euclidean}
A non-associative algebra $(A, \bullet)$ is called an \emph{anti-quadartic non-associative algebra} if there exists a pre-symplectic structure $\omega$ on $A$ that is anti-associative (or anti-invariant). In this case, the triple $(A, \bullet, \omega)$ denotes such a structure, and $\omega$ is called a \emph{anti-quadratic pre-symplectic form} on $(A, \cdot)$.
\end{Definition}

\sssbegin{Remark}
Let $(A, \bullet)$ be a non-associative algebra endowed with a non-degenerate \emph{skew-symmetric} (resp. \emph{symmetric}) bilinear form $\omega$. If $\omega$ is invariant (resp. anti-invariant), that is, $\omega(u \bullet v, w) = \omega(u, v \bullet w)$ for all $u, v, w \in A$, then $(A, \bullet)$ is necessarily trivial. 

We will give the prove in the case where the form is skew-symmetric. Indeed, for all $u, v, w \in A$, the combination of invariance and skew-symmetry yields:
\begin{align*}
\omega(u \bullet v, w) &= \omega(u, v \bullet w)= -\omega(v \bullet w, u)= -\omega(v, w \bullet u)  \\
&= \omega(w \bullet u, v) = \omega(w, u \bullet v)
= -\omega(u \bullet v, w).
\end{align*}
This implies that $2\omega(u \bullet v, w) = 0$, and since the characteristic of the base field is assumed to be different from 2, we obtain $\omega(u \bullet v, w) = 0$ for all $u, v, w \in A$. Finally, the non-degeneracy of $\omega$ ensures that $u \bullet v = 0$ for all $u, v \in A$, i.e., $(A, \bullet)$ is trivial.

It is worth noting that this situation contrasts sharply with the case where $\omega$ is symmetric. Nonassociative algebras endowed with a nondegenerate, invariant, and symmetric bilinear form (commonly referred to as \emph{quadratic algebras}) constitute a broad and extensively studied class of algebras, encompassing Lie, Jordan, associative, nearly associative, and Malcev algebras, among many others. \end{Remark}



In \cite{BE1}, the authors studied anti-quadratic Lie algebras and showed that such Lie algebras are $2$-step nilpotent.

\sssbegin{Proposition}
Let $(A, \bullet, \omega)$ be an anti-quadratic associative algebra. Then, $(A, \bullet)$ is nilpotent of index at most $3$.
\end{Proposition}

\begin{proof}
Since $\omega$ is anti-invariant, we have $\Ll_u^* = -\Rr_u$ for all $u \in A$. By associativity of $(A, \bullet)$, we have $\Ll_{u \bullet v} = \Ll_u \circ \Ll_v$ and $\Rr_{u \bullet v} = \Rr_v \circ \Rr_u$. Taking the adjoint yields
\[
\Ll_{u \bullet v}^* = (\Ll_u \circ \Ll_v)^* = \Ll_v^* \circ \Ll_u^*,
\]
which can be rewritten to 
\[
-\Rr_{u \bullet v}  = \Rr_v \circ \Rr_u = \Rr_{u \bullet v}.
\]
This implies $-\Rr_{u \bullet v} = \Rr_{u \bullet v}$, and thus $\Rr_{u \bullet v} = 0$ for all $u, v \in A$. Consequently, for any $w \in A$, we obtain:
\[
(w \bullet u) \bullet v  =w\bullet (u \bullet v)= 0.
\]
Hence, $A \bullet A \bullet A = \{0\}$, which means $(A, \bullet)$ is nilpotent of index $3$.
\end{proof}
\subsection{Nearly associative algebras}
Following \cite{DS}, a non-associative algebra $(A, \bullet)$ is said to be \emph{nearly associative} if, for all $u,v,w \in A$, the following identity holds:
\begin{equation}\label{nearly-associative}
(u\bullet v)\bullet w = v\bullet (w\bullet u).
\end{equation}
We refer to \cite[Example 3.17 and 3.20]{BBR} for examples of a nearly associative algebra that is neither commutative nor anti-commutative nor associative. 

If $(A, \bullet)$ is a nearly associative algebra, then for all $u,v \in A$, the left and right multiplication operators satisfy the following relations:
\begin{align}
\Ll_u \circ \Ll_v &= \Rr_v \circ \Rr_u, \label{left-right-1}\\
\Ll_u \circ \Rr_v &= \Ll_{v\bullet u}, \label{left-right-2}\\
\Rr_u \circ \Ll_v &= \Rr_{u\bullet v}. \label{left-right-3}
\end{align}
Conversely, if a non-associative algebra $(A, \bullet)$ satisfies any  of the conditions \eqref{left-right-1}, \eqref{left-right-2}, or \eqref{left-right-3}, then $(A, \bullet)$ is nearly associative.

\sssbegin{Lemma}[\cite{BBR}]
Let $(A, \bullet)$ be a nearly associative algebra. Then, for all $u,v, w,z \in A$, we have 
\begin{equation}\label{E1}
(u \bullet v) \bullet (w \bullet z) = (v \bullet u) \bullet (z \bullet w).
\end{equation}
\end{Lemma}

\sssbegin{Proposition}\cite[Propositions~4.2 and~4.5]{BBR}\label{Lie-Jordan}
Let $(A, \bullet)$ be a nearly associative algebra. Then the following statements hold:
\begin{enumerate}[$(i)$]
    \item  $A^{-}$ is a solvable Lie algebra.
    \item $A^{+}$ is a Jordan algebra.
\end{enumerate}
\end{Proposition}

In what follows, we introduce other classes of non-associative algebras that play a crucial role in the study of anti-quadratic nearly associative  algebras.

Following \cite{BE}, a non-associative algebra $(A, \bullet)$ is called an:
\begin{itemize}
    \item [$(i)$]\emph{anti-$\Rr$-algebra} if $\Rr_u\circ \Rr_v = -\Rr_v\circ \Rr_u $ for all $u, v \in A$, which is equivalent to 
    \begin{equation}\label{R-alg}
    (u\bullet v)\bullet w = -(u\bullet w)\bullet v , \quad \text{for all }\, u,v,w \in A.
    \end{equation}
    \item[$(ii)$] \emph{anti-$\Ll$-algebra} if $\Ll_u\circ \Ll_v  = -\Ll_v\circ \Ll_u $ for all $u,v \in A$, which means that
    \begin{equation}\label{L-alg}
    u\bullet(v\bullet w) = -v\bullet(u\bullet w), \quad \text{for all }\, u,v,w \in A.
    \end{equation}
    \item[$(iii)$] \emph{anti-$\Ll\Rr$-algebra} if $A$ is simultaneously an anti-$\Ll$-algebra and an anti-$\Rr$-algebra.
\end{itemize}

\sssbegin{Remark}
(i) Any associative anti-$\Ll\Rr$-algebra is a nearly associative algebra.

(ii) An algebra $(A, \bullet)$ satisfying the identities 
\[
\Rr_u\circ \Rr_v = \Rr_v\circ \Rr_u, \quad \Ll_u\circ \Ll_v  = \Ll_v\circ \Ll_u \quad \text{for all } u,v \in A, 
\]
is called bi-commutative; see, for example, \cite{IMS} and the references therein. 
\end{Remark}
\section{Structural Results on Anti-Quadratic Nearly Associative Algebras}\label{S3}

The purpose of this section is to characterize anti-quadratic nearly associative algebras and investigate their structural properties. We begin by establishing that, within the anti-quadratic framework, the classes of nearly associative algebras, anti-$\Ll$-algebras, and anti-$\Rr$-algebras completely coincide. Furthermore, we prove a striking structural result showing that any anti-quadratic nearly associative algebra is nilpotent of index at most $6$.  Building upon these algebraic features, we subsequently develop the representation theory for nearly associative algebras. Specifically, we define both the adjoint and coadjoint representations and determine the necessary and sufficient conditions for the existence of the latter. Finally, we establish that a nearly associative algebra admits an anti-quadratic bilinear form if and only if its adjoint and coadjoint representations are isomorphic $A$-bimodules.

\ssbegin{Proposition}
Let $(A, \bullet)$ be a nearly associative algebra. Then $(A, \bullet)$ is an anti-$\Rr$-algebra if and only if it is an anti-$\Ll$-algebra. 
\end{Proposition}

\begin{proof}
Suppose  that $(A, \bullet)$ is an anti-$\Rr$-algebra.  For any $u, v, w \in A$, we have:
\[
u \bullet (v \bullet w) = (w \bullet u) \bullet v = -(w \bullet v) \bullet u = -v \bullet (u \bullet w).
\]
This immediately yields $u \bullet (v \bullet w) = -v \bullet (u \bullet w)$, and hence $(A, \bullet)$ is an anti-$\Ll$-algebra. 

The converse follows with a similar argument.
\end{proof}

Now, we provide a characterization of anti-quadratic nearly associative algebras. 
\ssbegin{Theorem}\label{cara-nearly}
Let $(A, \bullet, \omega)$ be an anti-quadratic algebra. Then the following conditions are equivalent:
\begin{enumerate}
    \item[$(1)$] $(A, \bullet)$ is an anti-$\Ll$-algebra;
    \item[$(2)$] $(A, \bullet)$ is an anti-$\Rr$-algebra;
    \item[$(3)$] $(A, \bullet)$ is a nearly associative algebra.
\end{enumerate}
\end{Theorem}
\begin{proof}
We first prove that condition~(1) is equivalent to condition~(2). For all $u, v, w, z \in A$, using the properties of the anti-quadratic bilinear form $\omega$, we have the following: 
\begin{align*}
    \omega(u \bullet (v \bullet w) + v \bullet (u \bullet w), z)
    &= -\omega(u, (v \bullet w) \bullet z) - \omega(v, (u \bullet w) \bullet z) \\
    &= \omega((v \bullet w) \bullet z, u) + \omega((u \bullet w) \bullet z, v) \\
    &= \omega(v, w \bullet (z \bullet u)) + \omega(u, w \bullet (z \bullet v)) \\
    &= \omega(w, (z \bullet u) \bullet v) + \omega(w, (z \bullet v) \bullet u) \\
    &= -\omega((z \bullet u) \bullet v - (z \bullet v) \bullet u, w).
\end{align*}
By the non-degeneracy of $\omega$, it follows that $u \bullet (v \bullet w) + v \bullet (u \bullet w) = 0$ if and only if $(z \bullet u) \bullet v + (z \bullet v) \bullet u = 0$. This demonstrates that $(A, \bullet)$ is an anti-$\Ll$-algebra if and only if it is an anti-$\Rr$-algebra.

Next, we show that condition~(1) is equivalent to condition~(3). For all $u, v, w, z \in A$, we have:
\begin{align*}
    \omega(u \bullet (v \bullet w) + v \bullet (u \bullet w), z)
    &= -\omega(u, (v \bullet w) \bullet z) - \omega(v, (u \bullet w) \bullet z) \\
    &= -\omega(u, (v \bullet w) \bullet z) + \omega((u \bullet w) \bullet z, v) \\
    &=- \omega(u, (v \bullet w) \bullet z) + \omega(u, w \bullet (z \bullet v)) \\
    &= \omega((v \bullet w) \bullet z - w \bullet (z \bullet v), u).
\end{align*}
Again, by the non-degeneracy of $\omega$, this shows that condition~(1) is equivalent to condition~(3), which completes the proof.
\end{proof}

\ssbegin{Corollary}
Let $(A, \bullet, \omega)$ be an anti-quadratic nearly associative algebra. Then the following identities are satisfied:
\begin{equation}\label{cara1}
(u \bullet v) \bullet w = v \bullet (w \bullet u) = -w \bullet (v \bullet u) = -(u \bullet w) \bullet v, \quad \forall\, u, v, w \in A.
\end{equation}
\end{Corollary}

\begin{proof}
The result follows immediately from the equivalences established in Theorem \ref{cara-nearly}.
\end{proof}

\ssbegin{Theorem}\label{nil}
Let $(A, \bullet, \omega)$ be an anti-quadratic nearly associative algebra. Then, $(A, \bullet)$ is nilpotent of index at most $6$.
\end{Theorem}

\begin{proof}
Let $u, v, w, z, x, y \in A$. Using successively relations \eqref{E1}, \eqref{cara1}, and the definition of nearly associativity, we obtain
\begin{align*}
((u \bullet v) \bullet (w \bullet z)) \bullet (x \bullet y)
&\overset{\eqref{E1}}{=} ((w \bullet z) \bullet (u \bullet v)) \bullet (y \bullet x)
&\overset{\eqref{cara1}}{=} -((u \bullet v) \bullet (z \bullet w)) \bullet (y \bullet x) \\[4pt]
&\overset{\eqref{nearly-associative}}{=} - (z \bullet w) \bullet ((y \bullet x) \bullet (u \bullet v)) 
&\overset{\eqref{E1}}{=} - (z \bullet w) \bullet ((x \bullet y) \bullet (v \bullet u)) \\[4pt]
&\overset{\eqref{nearly-associative}}{=} - ((v \bullet u) \bullet (z \bullet w)) \bullet (x \bullet y) 
&\overset{\eqref{E1}}{=} - ((u \bullet v) \bullet (w \bullet z)) \bullet (x \bullet y).
\end{align*}
Thus,
\[
((u \bullet v) \bullet (w \bullet z)) \bullet (x \bullet y)
= - ((u \bullet v) \bullet (w \bullet z)) \bullet (x \bullet y),
\]
This chain of equalities implies that:
\[
2 \big( ((u \bullet v) \bullet (w \bullet z)) \bullet (x \bullet y) \big) = 0.
\]
Assuming the characteristic of the base field is different from $2$, it follows that $((u \bullet v) \bullet (w \bullet z)) \bullet (x \bullet y) = 0$. Consequently, any product of six elements in $A$ vanishes, which proves that $(A, \bullet)$ is nilpotent of index at most $6$.
\end{proof}

\ssbegin{Proposition}\label{Lie-Ja}
Let $(A, \bullet, \omega)$ be an anti-quadratic nearly associative algebra. Then:
\begin{enumerate}
    \item[$(1)$] $A^-$ is an anti-quadratic Lie algebra that is $2$-step nilpotent.
    \item[$(2)$] $A^+$ is an anti-quadratic Jacobi–Jordan algebra that is $2$-step nilpotent.
\end{enumerate}
\end{Proposition}

\begin{proof}
(1) By Proposition~\ref{Lie-Jordan}, $A^{-}$ is a Lie algebra. For all $u, v, w \in A$, the anti-quadratic property of $\omega$ yields
\[
\omega([u, v], w) 
= \omega(u \bullet v - v \bullet u, w) 
= -\omega(u, v \bullet w) + \omega(v, u \bullet w) 
= -\omega(u, v \bullet w - w \bullet v)
= -\omega(u, [v, w]).
\]
Hence, $(A^-, \omega)$ is anti-quadratic. Furthermore, using the definition of the adjoint operator and the Jacobi identity, we have $\mathrm{ad}_{[u, v]} = [\mathrm{ad}_u, \mathrm{ad}_v]$. Since $\omega$ is anti-quadratic, $\mathrm{ad}_u$ is symmetric, which implies $([\mathrm{ad}_u, \mathrm{ad}_v])^* = -[\mathrm{ad}_u, \mathrm{ad}_v]$ is skew-symmetric. Since $\ad_{[u,v]}$ is symmetric, we deduce that $\mathrm{ad}_{[u, v]} = 0$, meaning that $A^-$ is $2$-step nilpotent.

(2) According to Eqns~\eqref{cara1}, we have $(u \bullet v) \bullet w = -w \bullet (v \bullet u)$ for all $u, v, w \in A$. Expanding the product in $A^+$, we get:
\[
\{\{u, v\}, w\}
= (u \bullet v) \bullet w + (v \bullet u) \bullet w + w \bullet (u \bullet v) + w \bullet (v \bullet u)
= 0.
\]
This shows that $A^+$ is $2$-step nilpotent, which automatically satisfies the Jacobi-Jordan identity. The fact that $(A^+, \omega)$ is anti-quadratic follows from a similar computation as in the first part. This completes the proof.
\end{proof}

\ssbegin{Definition}[\cite{BBR1, DS}]\label{def-bimodule}
Let $(A, \bullet)$ be a nearly associative algebra, $V$ be a vector space, and $l, r: A \to \mathrm{End}(V)$ be two linear maps. The triple $(l, r, V)$ is called a \emph{bimodule} of $A$ if, for all $u, v \in A$, the following relations hold:
\begin{align}
l(u)l(v) &= r(v)r(u), \label{bim1} \\
l(u)r(v) &= l(v \bullet u), \label{bim2} \\
r(u)l(v) &= r(u \bullet v). \label{bim3}
\end{align}
\end{Definition}

\ssbegin{Example}
Let $(A, \bullet)$ be a nearly associative algebra. The triple $(A, \Ll, \Rr)$ is a bimodule of $A$, called the \emph{regular bimodule}, where $ \Ll, \Rr$ are left and right multiplications of $(\A, \bullet)$
\end{Example}

\ssbegin{Proposition}
Let $(l, r, V)$ be a bimodule of a nearly associative algebra $(A, \bullet)$. Then, the direct sum $A \oplus V$ can be endowed with a nearly associative algebra structure via the product defined by:
\begin{equation}\label{split-prod}
(u + x) \ast (v + y) := u \bullet v + l(u)y + r(v)x,
\end{equation}
for all $u, v \in A$ and $x, y \in V$.
\end{Proposition}

Let $A$ be an algebra and $(V, l, r)$ be an $A$-bimodule. As usual, we denote by $V^*$ the dual vector space of $V$. Then, the dual maps  $l^*: A \to \mathrm{End}(V^*)$  and $r^*: A \to \mathrm{End}(V^*)$ of $l$ and $r$, respectively, are defined for all $u \in A$ and $f \in V^*$ by:
\begin{equation}\label{dual-maps}
r^*(u)(f) := -f \circ l(u) \quad \text{and} \quad l^*(u)(f) := -f \circ r(u).
\end{equation}
Here, the composition of a linear form in $V^*$ with an endomorphism on $V$ is denoted by juxtaposition.

\ssbegin{Remark}
Let $(A, \bullet)$ be a nearly associative algebra and $(V, l, r)$ be an $A$-bimodule. If $r^*$ and $l^*$ are the dual maps of $r$ and $l$ defined by \eqref{dual-maps}, then $(V^*, r^*, l^*)$ is not necessarily an $A$-bimodule.
\end{Remark}

\ssbegin{Proposition}\label{Le1}
Let $(A, \bullet)$ be a nearly associative algebra and $(V, l, r)$ be its bimodule.. Then, the triple $(A^*, r^*, l^*)$ defines a bimodule \textup{(}representation\textup{)} of $A$ if and only if 
\begin{equation}\label{eq:rep-condition}
l(u \bullet v) = -r(v \bullet u), \quad \text{for all }\, u,v \in A.
\end{equation}
Under this condition, $(r^*, l^*, A^*)$ is called the \emph{anticoadjoint representation} \textup{(}or \emph{antico-regular representation}\textup{)} of $A$.
\end{Proposition}

\begin{proof}
Let $u, v \in A$ and $f \in A^*$. Using the definitions of the dual maps, we have
\[
\begin{array}{lcl}
r^*(u) \circ r^*(v)(f) 
&= & -r^*(u)(f \circ l(v)) 
= f \circ l(v) \circ l(u)
\overset{\eqref{bim1}}{=} f \circ r(u) \circ r(v) 
= -l^*(v)(f \circ r(u)) \\[2mm]
& = & l^*(v) \circ l^*(u)(f).
\end{array}
\]
This shows that the first identity $r^*(u) \circ r^*(v) = l^*(v) \circ l^*(u)$ always holds. Next, we have
\[
r^*(u) \circ l^*(v)(f) 
= -r^*(u)(f \circ r(v)) 
= f \circ r(v)\circ l(u) 
\overset{\eqref{bim3}}{=} f \circ r(v \bullet u).
\]
On the other hand, the definition gives $r^*(v \bullet u)(f) = -f \circ l(v \bullet u)$. Therefore, $r^*(u) \circ l^*(v)(f) = r^*(v \bullet u)(f)$ if and only if $f \circ r(v \bullet u) = -f \circ l(u \bullet v)$, which is equivalent to condition \eqref{eq:rep-condition}.

Similarly, for the last identity, we obtain 
\[
l^*(u) \circ r^*(v)(f) 
= -l^*(u)(f \circ l(v)) 
= f \circ l(v) \circ r(u) 
\overset{\eqref{bim2}}{=} f \circ l(u \bullet v),
\]
while $l^*(u \bullet v)(f) = -f \circ r(u \bullet v)$. Thus, $l^*(u) \circ r^*(v)(f) = l^*(u \bullet v)(f)$ if and only if condition \eqref{eq:rep-condition} holds. This completes the proof.
\end{proof}

\ssbegin{Proposition}\label{Isom}
Let $(A, \bullet)$ be a nearly associative algebra and $(\Ll, \Rr)$ denote its adjoint representation. Then the following conditions are equivalent:
\begin{enumerate}
    \item[$(1)$] $(\Ll^*, \Rr^*)$ is a representation of $A$ in $A^*$;
    \item[$(2)$] $(A, \bullet)$ is an anti-$\Ll$-algebra;
    \item[$(3)$]  $(A, \bullet)$ is an anti-$\Rr$-algebra.
\end{enumerate}
\end{Proposition}
\begin{proof}
Suppose first that $(\Ll^*, \Rr^*)$ is a representation of $A$ in $A^*$. Then, by Proposition~\ref{Le1}, we have $\Ll_{u \bullet v} = -\Rr_{v \bullet u}$ for all $u, v \in A$. Using the nearly associative identities, it follows that
\[
u \bullet (v \bullet w) = -(w \bullet v) \bullet u = -v \bullet (u \bullet w), \quad \forall\, u, v, w \in A,
\] 
which means that $A$ is an anti-$\Ll$-algebra.

Conversely, suppose that $A$ is an anti-$\Ll$-algebra. Then, we have
\[
(u \bullet v) \bullet w = v \bullet (w \bullet u) = -w \bullet (v \bullet u), \quad \forall\, u, v, w \in A,
\] 
which implies that $\Ll_{u \bullet v} = -\Rr_{v \bullet u}$ for all $u, v \in A$. Thus, by Proposition~\ref{Le1}, $(\Ll^*, \Rr^*)$ is a representation of $A$ in $A^*$. This proves the equivalence between (1) and (2).

Finally, the equivalence between conditions (2) and (3) follows directly from Proposition~\ref{cara-nearly}.
\end{proof}

\ssbegin{Corollary}
Let $(A, \bullet, \omega)$ be an anti-quadratic nearly associative algebra and $(\Ll, \Rr)$ denote its adjoint representation. Then, the pair $(\Ll^*, \Rr^*)$ defines a representation of $A$ in $A^*$.
\end{Corollary}

\begin{proof}
According to Theorem~\ref{cara-nearly}, $(A, \bullet)$ is both an anti-$\Ll$-algebra and an anti-$\Rr$-algebra. This structurally implies that
\[
\Ll_{u \bullet v} = -\Rr_{v \bullet u}, \quad \text{ for all }\, u, v \in A.
\]
Hence, by Proposition~\ref{Isom}, the pair $(\Rr^*,\Ll^*)$ indeed defines a representation of $A$ in $A^*$.
\end{proof}

\ssbegin{Proposition}
Let $(A, \bullet, \omega)$ be an anti-quadratic algebra, and let $\Ll$ and $\Rr$ denote the operators of left and right multiplication, respectively. Then, there exists an isomorphism $\Phi : A \rightarrow A^*$ such that
\[
\Phi(\Ll_u(a)) = \Ll^*_u(\Phi(a)) \quad \text{and} \quad \Phi(\Rr_u(a)) = \Rr^*_u(\Phi(a)), \quad \forall\, u, a \in A.
\]
\end{Proposition}

\begin{proof}
Consider the linear map $\Phi : A \rightarrow A^*$ defined for all $a \in A$ by $\Phi(a) := \omega(a, \cdot)$. Since $\omega$ is non-degenerate, $\Phi$ is an isomorphism of vector spaces. 

Moreover, for any $u, a, b \in A$, using the anti-associative of $\omega$ and the definitions of the dual maps, we have:
\[
\Phi(\Ll_u(a))(b) = \omega(u \bullet a, b) = -\omega(u, a \bullet b) = \omega(a \bullet b, u) = \omega(a, b \bullet u) = \Phi(a)(\Rr_u(b)) = \Ll^*_u(\Phi(a))(b).
\]
This immediately yields $\Phi(\Ll_u(a)) = \Ll^*_u(\Phi(a))$ for all $u, a \in A$. 

Similarly, an analogous computation shows that $\Phi(\Rr_u(a)) = \Rr^*_u(\Phi(a))$ for all $u, a \in A$, which completes the proof.
\end{proof}

\ssbegin{Theorem}\label{antieqdual}
Let $(A, \bullet)$ be nearly associative  algebra, and let $\Ll$ and $\Rr$ denote the operators of left and right multiplication, respectively. Then, there exists a bilinear form $\omega : A \times A \to \mathbb{K}$ for which $(A, \bullet, \omega)$ is anti-quadratic if and only if the dimension of $\A$ is even and there exists a vector space isomorphism $\Phi : A \to A^*$ satisfying:
\[
\Phi(\Ll_u(a)) = \Ll_u^*(\Phi(a)) \quad \text{and} \quad \Phi(\Rr_u(a)) = \Rr_u^*(\Phi(a)), \qquad \text{ for all }\, u,a \in A.
\]
\end{Theorem}
 \begin{proof}
If $(A, \bullet, \omega)$ is an anti-quadratic algebra, the conclusion follows immediately from Proposition~\ref{Isom}. Furthermore, since $\omega$ is a nondegenerate skew-symmetric bilinear form, $\dim A$ is necessarily even.  

Conversely, assume there exists a vector space isomorphism $\Phi : A \to A^*$ satisfying 
\[
\Phi(\Rr_u(a)) = \Rr_u^*(\Phi(a)) \quad \text{and} \quad \Phi(\Ll_u(a)) = \Ll_u^*(\Phi(a)), \quad \text{for all $u, a \in A$.}
\]
Define a bilinear form $M : A \times A \to \mathbb{K}$ by 
$
M(u,v) := \Phi(u)(v),$ for all $ u, v \in A.$ Since $\Phi$ is an isomorphism, then $M$ is non-degenerate. 

Moreover, for all $u, v, w \in A$, we have
\[
\begin{aligned}
M(u \bullet v, w) &= \Phi(u \bullet v)(w) = \Phi(\Rr_v(u))(w) = \Rr_v^*(\Phi(u))(w) \\
&= -\Phi(u)(\Ll_v(w)) = -\Phi(u)(v \bullet w) = -M(u, v \bullet w),
\end{aligned}
\]
which shows that $M$ is anti-invariant. Furthermore, $M$ satisfies the following identity:
\begin{equation}\label{eq:inv}
M(u \bullet v, w) = -M(v, w \bullet u), \qquad \forall\, u, v, w \in A.
\end{equation}
Indeed, for all $u, v, w \in A$, we have 
\[
\begin{aligned}
M(u \bullet v, w) &= \Phi(u \bullet v)(w) = \Phi(\Ll_u(v))(w) = \Ll_u^*(\Phi(v))(w) \\
&= -\Phi(v)(\Rr_u(w)) = -\Phi(v)(w \bullet u) = -M(v, w \bullet u).
\end{aligned}
\]
Let us now decompose $M$ into its symmetric and skew-symmetric parts, denoted by $M_s$ and $M_a$ respectively:
\[
M_s(u,v) := \frac{1}{2}\big( M(u,v) + M(v,u) \big), \qquad M_a(u,v) := \frac{1}{2}\big( M(u,v) - M(v,u) \big).
\]
Clearly, $M = M_s + M_a$. Since $M$ is anti-invariant and satisfies Eq.~\ref{eq:inv}, the symmetric part satisfies the following 
\[
\begin{aligned}
M_s(u \bullet v, w) &= \frac{1}{2}\big(M(u \bullet v, w) + M(w, u \bullet v)\big) \\
&= -\frac{1}{2}\big(M(u, v \bullet w) + M(v \bullet w, u)\big) = -M_s(u, v \bullet w).
\end{aligned}
\]
Hence, $M_s$ is a symmetric and an anti-invariant bilinear form. Analogously, $M_a$ is a skew-symmetric and an anti-invariant bilinear form.

Moreover, using both the symmetry and anti-invariance properties of $M_s$, we obtain, for all $u, v, w \in A$
\[
\begin{aligned}
M_s(u \bullet v, w) &= -M_s(u, v \bullet w) = -M_s(v \bullet w, u) = M_s(v, w \bullet u) \\ & =M_s(w \bullet u, v) = -M_s(w, u \bullet v) = -M_s(u \bullet v, w).
\end{aligned}
\]
Consequently, we obtain
\begin{equation}\label{eq:nul}
M_s(u \bullet v, w) = 0, \qquad \text{ for all }\, u, v, w \in A.
\end{equation}

Let $N$ denote the radical (kernel) of $M_s$, defined by
\[
N := \{u \in A \mid M_s(u,v) = 0, \, \forall\, v \in A\}.
\]
Equation~\ref{eq:nul} implies $A^2 \subseteq N$. Thus, $N$ is a two-sided ideal of $(A, \bullet)$.

Next, consider the radical of the skew-symmetric part:
\[
W := \{u \in A \mid M_a(u,v) = 0, \, \forall\, v \in A\}.
\]
The non-degeneracy of $M$ ensures that $N \cap W = \{0\}$. Therefore, there exists a vector subspace $V \subseteq A$ such that 
$
A = W \oplus V,$ with $N \subseteq V.$ Since $A^2 \subseteq N \subseteq V$, it follows that $V$ is also a two-sided ideal of $(A, \bullet)$, and the restriction $M_a|_{V \times V}$ is non-degenerate. Consequently, $V$ is even-dimensional. Since $A$ is also even-dimensional, $W$ is even-dimensional as well.

Let $H : W \times W \to \mathbb{K}$ be an arbitrary non-degenerate and skew-symmetric bilinear form. We extend $H$ to a bilinear form $\overline{H} : A \times A \to \mathbb{K}$ by setting
\[
\overline{H}|_{W \times W} = H \quad \text{and} \quad \overline{H}(V, A) = \overline{H}(A, V) = \{0\}.
\]
By construction, $\overline{H}$ is skew-symmetric. We now define the bilinear form $\omega$ on $A$ as follows
\[
\omega := M_a + \overline{H}.
\]
Since both $M_a$ and $\overline{H}$ are skew-symmetric, then $\omega$ is also skew-symmetric.

Let us prove that $\omega$ is non-degenerate. Let $w \in W$ and $v \in V$ be such that
$
\omega(w + v, u) = 0,$ for all $u \in A,$ that is equivalent to $\overline{H}(w, u) + M_a(v, u) = 0$.
\begin{itemize}
    \item If $u \in W$, then $M_a(v, u) = 0$. Hence, $\overline{H}(w, u) = 0$, which implies that $H(w, u) = 0$.  
Since $H$ is non-degenerate, it follows that $w = 0$.
    \item If $u \in V$, then $\overline{H}(w, u) = 0$. Hence, $M_a(v, u) = 0$.  
Since $M_a|_{V \times V}$ is non-degenerate, we obtain $v = 0$.  
\end{itemize}
Therefore, $\omega$ is non-degenerate.

Finally, since $A^2 \subseteq V$, it follows that $u \bullet v \in V$ and $v \bullet w \in V$ for all $u, v, w \in A$. Thus, $\overline{H}(u \bullet v, w) = \overline{H}(u, v \bullet w) = 0$. Using the anti-invariance property of $M_a$, we obtain
\[
\omega(u \bullet v, w) = M_a(u \bullet v, w) = -M_a(u, v \bullet w) = -\omega(u, v \bullet w).
\]
Therefore, $\omega$ is anti-invariant, and hence $(A, \bullet, \omega)$ is an anti-quadratic algebra.
\end{proof}

\section{\texorpdfstring{$TT^*$}{TT*-} extensions for nearly associative algebras}\label{S4}
We now introduce a novel construction inspired by the classical $TT^*$-extension method, originally developed by Bordemann \cite{Bordemann} for quadratic non-associative algebras. This method was later adapted to anti-quadratic Lie algebras in \cite{BE1} and is referred to as a Lagrangian extension in \cite{BC}; see also \cite{BM, BR} for the corresponding notion in the super setting. Our goal is to adapt the core philosophy of this extension to the setting of anti-quadratic nearly associative algebras, taking into account the structural peculiarities of the dual maps and the underlying anti-quadratic bilinear forms. We show that every anti-quadratic nearly associative algebra can be obtained via a $TT^*$-extension.

\ssbegin{Definition}[2-cocycle condition]\label{def-cocycle}
Let $(A, \bullet)$ be a nearly associative algebra, $(V, r, l)$ be an $A$-bimodule, and $\theta : A \times A \to V$ be a bilinear map. Then, $\theta$ is called a $2$-cocycle on $A$ with values in $V$ if for all $u, v, w \in A$, the following identity holds
\begin{equation}\label{2-cocycle-cond}
\theta(u \bullet v, w) + r(w)(\theta(u,v)) = \theta(v, w \bullet u) + l(v)(\theta(w,u)).
\end{equation}
\end{Definition}
\ssbegin{Theorem}\label{TT-extension}
Let $(A, \bullet)$ be a nearly associative anti-$\Ll$-algebra and $\theta : A \times A \to A^*$ be a bilinear map. The vector space $TT^*_\theta (A) = A \oplus A^*$ endowed with the product defined by: 
\begin{equation}\label{pr-TT}
(u + f) \diamond  (v + g) := u \bullet v + \theta (u,v) + \Rr^*_v(f) + \Ll^*_u(g), \quad \text{ for all }\, u,v \in A, \; f, g \in A^*,
\end{equation}
is nearly associative if and only if $\theta $ is a $2$-cocycle on $A$ with values in the  representation $(\Ll^*, \Rr^*)$.

Moreover, the bilinear form defined by:
\[
\omega : TT^*_\theta (A) \times TT^*_\theta (A) \to \mathbb{K}, \quad \omega(u + f, v + g) = f(v) - g(u),
\]
is skew-symmetric and non-degenerate on $TT^*_\theta (A)$. Furthermore, $\omega$ is anti-associative on $(TT^*_\theta (A), \diamond )$ if and only if $\theta $ satisfies the following cyclic identity:
\begin{equation}
\theta (u,v)(w) = \theta (v,w)(u), \quad \text{ for all }\, u,v,w \in A.
\end{equation}
In this case, the anti-quadratic nearly associative algebra $(TT^*_\theta (A), \diamond , \omega)$ is called the $TT^*_\theta $-extension of $A$ by means of $\theta$.
\end{Theorem}

\begin{proof}
For every $u,v,w \in A$ and $f,g,h \in A^*$, we have
\begin{align*}
((u+f) \diamond (v+g)) \diamond (w+h) &= (u \bullet v) \bullet w + \theta(u \bullet v, w) + \Rr^*_{w}\big(\theta(u,v) + \Rr^*_v(f) + \Ll^*_u(g)\big)\\
& \; \;  + \Ll^*_{u \bullet v}(h), \\
(v+g) \diamond ((w+h) \diamond (u+f)) &= v \bullet (w \bullet u) + \theta(v, w \bullet u) + \Ll^*_{v}\big(\theta(w,u) + \Rr^*_u(h) + \Ll^*_w(f)\big)\\
& \; \;  + \Rr^*_{w \bullet u}(g).
\end{align*}
Since $(A, \bullet)$ is a nearly associative algebra and $(\Ll^*, \Rr^*)$ is a representation of $A$ on $A^*$, it follows that the product $\diamond$ is nearly associative if and only if
\[
\theta(u \bullet v, w) + \Rr^*_{w}(\theta(u,v)) = \theta(v, w \bullet u) + \Ll^*_{v}(\theta(w,u)).
\]
This means that the product is nearly associative if and only if $\theta$ is a $2$-cocycle on $A$ with values in the representation $(\Ll^*, \Rr^*)$.

\medskip
Moreover, by construction, $\omega$ is clearly a non-degenerate skew-symmetric bilinear form. Now, for all $u,v,w \in A$ and $f,g,h \in A^*$, we have 
\begin{align*}
\omega((u+f) \diamond (v+g), w+h) &= -h(u \bullet v) + \theta(u,v)(w) + \Rr^*_v(f)(w) + \Ll^*_u(g)(w) \\
&= -h(u \bullet v) + \theta(u,v)(w) - f(v \bullet w) - g(w \bullet u),
\end{align*}
and
\begin{align*}
\omega(u+f, (v+g) \diamond (w+h)) 
&= f(v \bullet w) - \theta(v,w)(u) - \Rr^*_w(g)(u) - \Ll^*_v(h)(u) \\
&= f(v \bullet w) - \theta(v,w)(u) + g(w \bullet u) + h(u \bullet v).
\end{align*}
It follows that $\omega$ is anti-invariant  if and only if
\[
\theta(u,v)(w) = \theta(v,w)(u).
\]
The proof is now complete.
\end{proof}

\ssbegin{Theorem}\label{TT-ideal}
Let $(A, \bullet, \omega)$ be an anti-quadratic nearly associative algebra of dimension $2n$. Then, $(A, \bullet, \omega)$ is isometrically isomorphic to a $TT^*$-extension $(TT^*_\Omega(B), \diamond , \omega_c)$ if and only if $A$ contains a Lagrangian ideal $I$. In this case, $B \cong A/I$.
\end{Theorem}
 Note that any Lagrangian subspace $I$ of $A$ is an ideal of $A$ if and only if $I \bullet I = 0$.
\begin{proof}
Assume that $A$ is $TT^*$-extension of $B$ by means of $\theta$. It is clear that $B^*$ is a Lagrangian subspace of $TT^*_\theta(B)$, and by the product \eqref{pr-TT}, it is easy to see that $B^*$ is a two-sided of $TT^*_\theta(B)$.

Now, assume that  $A$ contains a Lagrangian two-sided ideal $I$, i.e., $I=I^\bot$. If $I$ is a two-sided ideal of $A$, then by the anti-invariant of $\omega$ , we have $\omega(I \bullet I, A) = -\omega(I, I \bullet A)=0$, as $I$ is two-sided ideal and Lagrangian. Since $\omega$ is non-degenerate, it follows that $I \bullet I = 0$.  Conversely, if $I \bullet I = 0$, then $\omega(I \bullet A, I) = -\omega(I, A \bullet I)=\omega(A, I\bullet I)= 0$, which implies $ I\bullet A, A\bullet I \subseteq I^\perp = I$, showing that $I$ is a two-sided ideal of $A$.

Suppose now that $I$ is a Lagrangian two-sided ideal of $A$.  Let $B$ denote the  algebra $A/I$ and $p : A \to B$ be the canonical projection. We can choose a Lagrangian vector subspace $B_0$ to $I$ in $A$, i.e., $A = B_0 \oplus I$ and $B_0$ is Lagrangian. Let $p_0 : A \to B_0$ and $p_1 : A \to I$ denote the projections along $I$ and $B_0$, respectively. 

We define the linear map $\Phi : I \to B^*$ by $\Phi(u)(p(a)) := \omega(u, a)$. This is well-defined because $I$ is Lagrangian. Since $\omega$ is non-degenerate, $I^\perp = I$, and $\dim I  = \dim B$, the map $\Phi$ is a vector space isomorphism. Furthermore, $\Phi$ intertwines the algebra actions. Indeed, for any $a, x \in A$ and $u \in I$, the anti-quadratic property yields
\begin{align*}
\Phi(a \bullet u)(p(x)) &= \omega(a \bullet u, x) = -\omega(a, u \bullet x) = \omega(u\bullet x, a) =- \omega(u, x \bullet a) \\
&=- \Phi(u)(p(x) \bullet p(a)) = -\Phi(u)(\Rr_{p(a)}(p(x))) = \Ll^*_{p(a)}(\Phi(u))(p(x)).
\end{align*}
Similarly, 
$$
\Phi(u \bullet a)(p(x)) = \omega(u \bullet a, x) = -\omega(u, a \bullet x) = -\Phi(u)(\Ll_{p(a)}(p(x))) = \Rr^*_{p(a)}(\Phi(u))(p(x)).
$$
Hence, $\Phi(a \bullet u) = \Ll^*_{p(a)}(\Phi(u))$ and $\Phi(u \bullet a) = \Rr^*_{p(a)}(\Phi(u))$.

Next, we define the bilinear map $\theta : B \times B \to B^*$ by
\[
\theta (p(b_0), p(b'_0)) := \Phi(p_1(b_0 \bullet b'_0)), \quad \text{for } b_0, b'_0 \in B_0.
\]
This is well-defined since the restriction of $p$ to $B_0$ is a linear isomorphism. Now, define the linear map $\varphi : A \to B \oplus B^*$ by $\varphi(b_0 + u) := p(b_0) + \Phi(u)$. Clearly, $\varphi$ is a linear isomorphism. We show that $\varphi$ is an algebra isomorphism from $(A, \bullet)$ to $(TT^*_\theta(B), \diamond)$. For $b_0, b'_0 \in B_0$ and $u, u' \in I$, using $I \bullet I = 0$, we have 
\begin{align*}
\varphi((b_0 + u) \bullet (b'_0 + u')) &= \varphi(p_0(b_0 \bullet b'_0) + p_1(b_0 \bullet b'_0) + b_0 \bullet u' + u \bullet b'_0) \\
&= p(b_0 \bullet b'_0) + \Phi(p_1(b_0 \bullet b'_0)) + \Phi(b_0 \bullet u') + \Phi(u \bullet b'_0) \\
&= p(b_0) \bullet p(b'_0) + \theta (p(b_0), p(b'_0)) + \Ll^*_{p(b_0)}(\Phi(u')) + \Rr^*_{p(b'_0)}(\Phi(u)) \\
&\overset{\eqref{pr-TT}}{=} (p(b_0) + \Phi(u)) \diamond (p(b'_0) + \Phi(u')) \\
&= \varphi(b_0 + u) \diamond \varphi(b'_0 + u').
\end{align*}

Finally, we verify that $\varphi$ preserves the anti-quadratic forms. Since $B_0$ and $I$ are both Lagrangian, we get:
\begin{align*}
(\varphi^*\omega_c)(b_0 + u, b'_0 + u') &= \omega_e(p(b_0) + \Phi(u), p(b'_0) + \Phi(u')) \\
&= \Phi(u)(p(b'_0)) - \Phi(u')(p(b_0)) \\
&= \omega(u, b'_0) - \omega(u', b_0) \\
&= \omega(b_0, b'_0) + \omega(b_0, u') + \omega(u, b'_0) + \omega(u, u') \\
&= \omega(b_0 + u, b'_0 + u').
\end{align*}
Hence, $\varphi^*\omega_c = \omega$. This implies that $\omega_c$ is an anti-quadratic form on $TT^*_\theta (B)$, which in turn forces $\theta $ to be a cyclic $2$-cocycle. Thus, $(A, \bullet, \omega)$ and $(TT^*_\Omega(B), \diamond, \omega_c)$ are isometrically isomorphic as anti-quadratic nearly associative algebras.
\end{proof}

\ssbegin{Proposition}\label{Propriete}
Let $(A, \bullet, \omega)$ be an anti-quadratic nearly associative algebra. Then the following assertions hold:
\begin{enumerate}
    \item[$(1)$] $N(A) = (A \bullet A)^\perp$,
    \item[$(2)$] If $I$ is a two-sided ideal of $(A, \bullet)$, then its orthogonal complement $I^\perp$ is also a two-sided ideal of $(A, \bullet)$. Moreover, they satisfy the properties:
    $
    I^\perp \bullet I = I \bullet I^\perp = \{0\}$,

    \item[$(3)$]  $N(A)$ and $A \bullet A$ are both degenerate subspaces of $A$ with respect to $\omega$. 
    \item[$(4)$]  The space $I := N(A) \cap N(A)^\perp$ and its orthogonal complement $I^\perp$ are both two-sided ideals of  $(A, \bullet)$.

\end{enumerate}
\end{Proposition}

\begin{proof}
Let us proof Part (1). Let $u \in N(A)$. For all $v,w  \in A$, and by the anti-invariance of $\omega$, we have
    $$\omega(u, v \bullet w) = -\omega(u \bullet v, w) = 0.$$ Hence $u \in (A \bullet A)^\perp$. Conversely, suppose  $u \in (A \bullet A)^\perp$. For all $v,w\in A$, we have  $$\omega(v\bullet u, w)=-\omega(v, u\bullet w)=-\omega(w\bullet v, u)=0.$$ By the non-degeneracy of $\omega$, it follows that $v \bullet u = 0$ and $u\bullet w=0$, for all $v,w\in A$. Hence, $u \in N(A)$.
    
Let us proof Part (2). Let $v \in I^\perp$, $u \in A$, and $i \in I$. We have 
$$
\omega(u \bullet v, i) = -\omega(u, v \bullet i)=-\omega(v, i\bullet u)=0,\esp \omega(v\bullet u, i)=-\omega(v, u\bullet i)=0.
$$ 
It follows that $u \bullet v, v\bullet u \in I^\perp$. The fact that $I^\perp \bullet I = \{0\}$ follows directly because 
$$
\omega(I^\perp \bullet I, A) = -\omega(I^\perp, I \bullet A) = \{0\}.
$$
    
 If $I$ is a Lagrangian two-sided ideal of $(A, \bullet)$, i,e., $I = I^\perp$, from item 2 we have $I \bullet I = I^\perp \bullet I = \{0\}$. Conversely, if $I$ is Lagrangian  and $I \bullet I = \{0\}$, then for any $u \in A$, we have $\omega(u \bullet I, I) = -\omega(u, I \bullet I)=\{0\}$ and $\omega(I\bullet u, I) = -\omega(I, u \bullet I)=-\omega( u, I\bullet I)=\{0\}.$ Thus, we get $A \bullet I, I\bullet A \subseteq I$ and hence $I$ is two-sided ideal of $(A, \bullet)$. 

Let us prove Part (3). According to Theorem \ref{nil},  $(A, \bullet)$ is  nilpotent,  and therefore $N(A) \cap (A \bullet A) \neq \{0\}$. From Part (1), we established that $(A \bullet A)^\perp = N(A)$, and therefore $N(A) \cap N(A)^\perp \neq \{0\}$ and $(A \bullet A)^\perp \cap (A \bullet A) \neq \{0\}$. Consequently, both $N(A)$ and $A \bullet A$ are degenerate subspaces of $A$ with respect to $\omega$.
    
Let us prove Part (4). Since $N(A)$ is a two-sided ideal of $(A, \bullet)$, it follows from Part (2) that its orthogonal complement $N(A)^\perp$ is also a two-sided ideal. Because the intersection of two two-sided ideals is always a two-sided ideal, the subspace $I := N(A) \cap N(A)^\perp$ is a two-sided ideal. Furthermore, by Part (2)  its orthogonal complement $I^\perp$ is also a two-sided ideal of $(A, \bullet)$.
\end{proof}

\ssbegin{Proposition}\label{pr0}
Let $(A, \bullet, \omega)$ be an anti-quadratic nearly associative algebra. Then there exists a two-sided Lagrangian ideal $W_{\max}$ of $(A, \bullet)$.
\end{Proposition}

\begin{proof}
We proceed by induction on $n$, where $\dim A = 2n$. The base case $n = 0$ is trivially satisfied. Let us assume that the result holds for any anti-quadratic algebra of dimension strictly less than $2n$, and let $(A, \bullet, \omega)$ be a $2n$-dimensional anti-quadratic algebra.

If $(A, \bullet)$ is the trivial algebra (i.e., $A \bullet A = \{0\}$), the result is straightforward.

Now, assume that $(A, \bullet)$ is a non-trivial algebra. According to Theorem \ref{nil}, the algebra $(A, \bullet)$ is nilpotent, which implies that $N(A) \cap (A \bullet A) \neq \{0\}$. Proposition \ref{Propriete} implies that $N(A) = (A \bullet A)^\perp$. Consequently, $N(A) \cap N(A)^\perp \neq \{0\}$. Let us define the subspace $W := N(A) \cap N(A)^\perp$. It is straightforward to see that $W$ is totally isotropic, and for any $u \in W$, we have $\Ll_u = \Rr_u = 0$, which ensures that $W$ is a two-sided ideal of $(A, \bullet)$. 

Furthermore, for all $u \in W$, $v \in W^\perp$, and $z \in A$, the anti-invariance of $\omega$ yields
\[
\omega(z \bullet v, u) = -\omega(v, u \bullet z) = 0 \quad \text{and} \quad \omega(v \bullet z, u) = -\omega(v, z \bullet u) = 0.
\]
Thus, $W^\perp$ is also a two-sided ideal of $(A, \bullet)$.

Next, consider the quotient space $A_0 := W^\perp / W$, which is well-defined since $W$ is totally isotropic and $W \subseteq W^\perp$. Let $p: W^\perp \to A_0$ denote the canonical projection. The symplectic form $\omega$ and the product $\bullet$ naturally induce a well-defined symplectic form $\omega_0$ and a product $\bullet_0$ on $A_0$, respectively, defined by
\[
\omega_0(p(u), p(v)) := \omega(u,v) \quad \text{and} \quad p(u) \bullet_0 p(v) := p(u \bullet v) \quad \text{for all } u, v \in W^\perp.
\]
Thus, $(A_0, \bullet_0, \omega_0)$ inherits the structure of a nilpotent anti-quadratic nearly associative algebra. The dimension of this lower-dimensional quotient space is given by
\[
\dim A_0 = \dim W^\perp - \dim W = (2n - \dim W) - \dim W = 2(n - \dim W).
\]
By the induction hypothesis, there exists a Lagrangian two-sided ideal $W_0$ of $(A_0, \bullet_0)$, which satisfies
\[
\dim W_0 = \frac{1}{2}\dim A_0 = n - \dim W.
\]

We now define $W_{\max} := p^{-1}(W_0)$ as the preimage of $W_0$ under the canonical projection $p$. By construction, $W_{\max} \subseteq W^\perp$ contains $W$, and we have the vector space isomorphism $W_{\max}/W \cong W_0$. Consequently, the dimension of $W_{\max}$ is given by
\[
\dim W_{\max} = \dim W_0 + \dim W = (n - \dim W) + \dim W = n.
\]
Let us show that $W_{\max}$ is totally isotropic. Let $u, v \in W_{\max}$. Since $p(u), p(v) \in W_0$ and $W_0$ is isotropic in $A_0$, we have
\[
\omega(u, v) = \omega_0(p(u), p(v)) = 0.
\]
 Given that $\dim W_{\max} = n$, it constitutes a Lagrangian subspace of $A$.

Finally, let us verify that $W_{\max}$ is a two-sided ideal of $(A, \bullet)$. For any $u \in W_{\max} \subseteq W^\perp$ and $v \in A$, since $W^\perp$ is a two-sided ideal, we have $u \bullet v \in W^\perp$ and $v \bullet u \in W^\perp$. Taking their projections, we obtain
\[
p(u \bullet v) = p(u) \bullet_0 p(v) \in W_0 \quad \text{and} \quad p(v \bullet u) = p(v) \bullet_0 p(u) \in W_0,
\]
because $p(u) \in W_0$ and $W_0$ is a two-sided ideal of $A_0$. This implies that $u \bullet v \in p^{-1}(W_0) = W_{\max}$ and $v \bullet u \in p^{-1}(W_0) = W_{\max}$. 

Therefore, $W_{\max}$ is a two-sided ideal of $(A, \bullet)$, and hence it is a Lagrangian two-sided ideal of $(A, \bullet, \omega)$, which completes the proof.
\end{proof}

The following theorem gives the converse of Theorem \ref{TT-extension}.
\ssbegin{Theorem}\label{Converse-T}
Let $(A, \bullet, \omega)$ be an anti-quadratic nearly associative algebra. Then $(A, \bullet, \omega)$ is a $TT^*_\theta$-extension of a nearly associative anti-$\Ll$-algebra $(B, \bullet_B)$ by means of $\theta$.
\end{Theorem}
\begin{proof}
By virtue of Proposition \ref{pr0}, the anti-quadratic nearly associative algebra $(\mathcal{A}, \bullet, \omega)$ admits a Lagrangian two-sided ideal $W_{\max}$. Consequently, applying Theorem \ref{TT-ideal}, we conclude that $(\mathcal{A}, \bullet, \omega)$ is a $T^*$-extension of $W_{\max}$ by means of $\theta$, which completes the proof.
\end{proof}

\section{Double extension of anti-quadratic nearly associative algebras}\label{S5}
In this section, we introduce the notions of annihilator extensions and almost semidirect products of nearly associative algebras. These two structural notion play a crucial role in our framework, as they constitute the foundational building blocks required to define and construct the notion of a \emph{double extension} for anti-quadratic nearly associative algebras. Furthermore, we prove that every anti-quadratic nearly associative algebra can be obtained by means of such a double extension.

\subsection{Annihilator extensions of nearly associative algebras}

First, we investigate the notion of an annihilator extension of a nearly associative algebra. let $H := \K e$ be a one-dimensional vector space, and let $\varphi : A \times A \longrightarrow \K$ be a bilinear form.  We define a new bilinear product $\Tilde{\bullet}$ by
\[
(u + \al e) \Tilde{\bullet} (v + \beta e) := u \bullet v + \varphi(u,v) e, \qquad \text{ for all }\, u,v \in A,\; \al, \beta \in \mathbb{K}.
\]
Then $(\widetilde{A}, \bullet)$ is a nearly associative algebra if and only if the bilinear map $\varphi$ satisfies the following condition:
\begin{equation}\label{line-central}
\varphi(u \bullet v, w) = \varphi(v, w \bullet u), \qquad \text{ for all}\, u,v,w \in A.
\end{equation}

If a bilinear map $\varphi$ satisfies Condition~\eqref{line-central}, then the nearly associative algebra $(\widetilde{A}, \bullet)$ is called the \emph{annihilator extension} of the nearly associative algebra $(A, \bullet)$ by means of $\varphi$.



\subsection{Almost semi-direct products of nearly associative algebras}

Lastly, we introduce the notion of an almost semi-direct product of nearly associative algebras. Let $(A, \bullet)$ be a nearly associative algebra,  let $H := \K d$ be a one-dimensional vector space. Let $\xi, D : A \longrightarrow A$  be two linear maps linear map,  let $a_0 \in A$ and a scalar $\la\in\mathbb{K}$. We define a bilinear product $\overline{\bullet}$ by:
\[
u \overline{\bullet} v = u \bullet v, \qquad
u \overline{\bullet} d = \xi(u), \qquad
d \overline{\bullet} u = D(u), \qquad
d \overline{\bullet} d = \lambda d + a_0,
\]
for all $u,v \in A$.

Then, $(\overline{A}, \overline{\bullet})$ is a nearly associative algebra if and only if the following compatibility conditions hold:
\begin{equation}\label{line-semi}
\begin{aligned}
\xi(u\bullet v)
&=v\bullet D(u),
&
\xi(u)\bullet v
&=D(v\bullet u),
\\[0.3em]
u\bullet \xi(v)
&=D(u)\bullet v,
&
\xi^2
&=D^2,
\\[0.3em]
\xi\circ D
&=\Rr_{a_0}+\lambda\,\xi,
&
D\circ\xi
&=\Ll_{a_0}+\lambda\,D,
\\[0.3em]
(\xi-D)(a_0)
&=0.
\end{aligned}
\end{equation}
for all $u,v \in A$ and $X,Y,Z \in V$.

If Conditions~\eqref{line-semi} are satisfied, the nearly associative algebra $(\overline{A}, \overline{\bullet})$ is called the \textit{almost semi-direct product} of the nearly associative algebra $(A, \bullet)$  by means of $( D, \xi,  a_0, \la)$. In this context, the tuple $(\mathbb{K} d, A, D, \xi, a_0, \la)$ is called an \textit{almost semi-direct product context} of nearly associative algebras.


%


\subsection{Double extension of anti-quadratic nearly associative algebras}

Let $(A, \bullet_A, \omega_A)$ be an anti-quadratic nearly associative algebra. Let $H := \K d$ be a one-dimensional vector space and let $H^* := \K e$ denote its dual space.  

In order to introduce the process of double extension by $H$, we consider a linear map $D: A \rightarrow A$, an element $a_0 \in A$, and a scalar $\mu \in \K$. We extend the linear map $D$ and the element $a_0$ to the linear maps $\widetilde{\xi}, \widetilde{D}: A \oplus H^* \rightarrow A \oplus H^*$ as follows:    
\begin{equation}\label{eq:double-tilde}
\widetilde{\xi}(u + \alpha e) := -D^*(u) + \omega_A(a_0, u)e,\;\;  \widetilde{D}(u + \alpha e) := D(u) + \omega_A(a_0, u)e, \;\; \widetilde{a}_0 = a_0 + \mu e,
\end{equation}
for all $u \in A$ and $\alpha \in \K$.

\medskip

\sssbegin{Lemma} \label{Lemma1}
The vector space $\widetilde{A} := A \oplus H^*$ equipped with the bilinear product
\[
(u + \alpha e) \widetilde{\bullet} (v + \beta e) := u \bullet_A v + \omega_A(D(u), v) e, \qquad \text{for all } u, v \in A, \; \alpha, \beta \in \K,
\]
is a nearly associative algebra if and only if
\begin{equation}\label{Le-centrale}
D(u \bullet_A v) = -u \bullet_A D(v), \qquad \text{for all } u, v \in A.
\end{equation}
\end{Lemma} 

\begin{proof}
Set $\phi(u,v) := \omega_A(D(u), v)$ for all $u, v \in A$. For any $u, v, w \in A$, we have:
\[
\phi(u \bullet_A v, w) - \phi(v, w \bullet_A u) = \omega_A(D(u \bullet_A v), w) - \omega_A(D(v), w \bullet_A u).
\]
Using the anti-invariance of $\omega_A$, the second term can be rewritten as:
\[
-\omega_A(D(v), w \bullet_A u) = \omega_A(w, u \bullet_A D(v)).
\]
It follows that 
\[
\phi(u \bullet_A v, w) - \phi(v, w \bullet_A u) = \omega_A\big(D(u \bullet_A v) + u \bullet_A D(v), w\big).
\]
Since $\omega_A$ is non-degenerate, it follows that $\phi$ satisfies Condition~\eqref{line-central} if and only if $D(u \bullet_A v) = -u \bullet_A D(v)$ for all $u, v \in A$, which completes the proof.
\end{proof}

\sssbegin{Lemma}\label{Lemma2}
$(H, \widetilde{A} := A \oplus H^*, \widetilde{D}, \widetilde{\xi}, \widetilde{a}_0, 0)$ is a context of almost semi-direct product of nearly associative algebras if and only if, for all $u, v = A$, the following conditions hold:
\begin{equation}\label{con-ex}
\begin{aligned}
D^*(u\bullet_A v)
&=-v\bullet_A D(u),
&
D^*(u)\bullet_A v
&=-D(v\bullet_A u),
\\[0.3em]
u\bullet_A D^*(v)
&=-D(u)\bullet_A v,
&
D^2
&=0,
\\[0.3em]
D^*\circ D
&=-D\circ D^*
=\Ll_{a_0},
&
D(a_0)
&=D^*(a_0)=0.
\end{aligned}
\end{equation}
\end{Lemma}

\begin{proof}
Let $u, v \in A$ and $\alpha, \beta \in \K$. A straightforward computation shows that the following equivalences hold in $(\widetilde{A},\widetilde{\bullet})$:
\begin{itemize}
    \item[$-$] The condition $\widetilde{\xi}((u+\alpha e)\widetilde{\bullet} (v+\beta e)) = (v+\beta e) \widetilde{\bullet} \widetilde{D}(u+\alpha e)$ holds if and only if
    \[
    D^*(u \bullet_A v) = -v \bullet_A D(u) \quad \text{and} \quad D^* \circ D = \Ll_{a_0}.
    \]

    \item[$-$] The equation $\widetilde{\xi}(u+\alpha e)\widetilde{\bullet} (v+\beta e) = \widetilde{D}((v+\beta e)\widetilde{\bullet} (u+\alpha e))$ is equivalent to
    \[
    D^*(u) \bullet_A v = -D(v \bullet_A u) \quad \text{and} \quad D \circ D^* = \Rr_{a_0}.
    \]

    \item[$-$] The relation $(u+\alpha e)\widetilde{\bullet} \widetilde{\xi}(v+\beta e) = \widetilde{D}(u+\alpha e) \widetilde{\bullet} (v+\beta e)$ is equivalent to
    \[
    u \bullet_A D^*(v) = -D(u) \bullet_A v \quad \text{and} \quad D^2 = 0.
    \]

    \item[$-$] The identity $\widetilde{\xi}^2(u+\alpha e) = \widetilde{D}^2(u+\alpha e)$ is equivalent to
    \[
    D^{*2} = D^2 \quad \text{and} \quad (D + D^*)(a_0) = 0.
    \]

    \item[$-$] The equation $\widetilde{\xi} \circ \widetilde{D}(u+\alpha e) = (u+\alpha e) \widetilde{\bullet} \widetilde{a}_0$ is equivalent to
    \[
    D^* \circ D = -\Rr_{a_0} \quad \text{and} \quad D^*(a_0) = 0.
    \]

    \item[$-$] The relation $\widetilde{D} \circ \widetilde{\xi}(u+\alpha e) = \widetilde{a}_0 \widetilde{\bullet} (u+\alpha e)$ is equivalent to
    \[
    D \circ D^* = -\Ll_{a_0} \quad \text{and} \quad D(a_0) = 0.
    \]

    \item[$-$] Finally, the equation $(\widetilde{D} - \widetilde{\xi})(\widetilde{a}_0) = 0$ is equivalent to
    \[
    (D + D^*)(a_0) = 0.
    \]
\end{itemize}
From these relations, we observe that $D^* \circ D = \Ll_{a_0}$ and $D \circ D^* = -\Ll_{a_0}$. Since $\omega_A$ is anti-invariant on $(A, \bullet_A)$, the left multiplication operator satisfies $\Ll_{a_0}^* = -\Rr_{a_0}$. Consequently, we deduce $D^* \circ D = -\Rr_{a_0}$ and $D \circ D^* = \Rr_{a_0}$. This completes the proof of the lemma.
\end{proof}

We next introduce the notion of a double extension for anti-quadratic nearly associative algebras. This construction, inspired by the classical double extension method of A. Medina and P. Revoy~\cite{Medina-Revoy} for quadratic Lie algebras, provides an effective way to generate new anti-quadratic structures from lower-dimensional ones. The precise formulation of this extension process is given in the following theorem.

\sssbegin{Theorem}\label{construction-ex}
Let $(A, \bullet_A, \omega_A) $ be an anti-quadratic nearly associative algebra, let $H = \K d$ be a one-dimensional vector space, and let $H^* = \K e$ be its dual space. Consider an endomorphism $D \in \mathrm{End}(A)$, an element $a_0 \in A$, and a scalar $\mu \in \K$ such that the following conditions are satisfied for all $u, v \in A$:
\begin{equation}\label{con-ex2}
\begin{aligned}
D(u\bullet_A v)
&=-u\bullet_A D(v), &
D^*(u\bullet_A v)
&=-v\bullet_A D(u),\\[0.3em]
D^*(u)\bullet_A v
&=-D(v\bullet_A u), &
u\bullet_A D^*(v)
&=-D(u)\bullet_A v,\\[0.3em]
D^2&=0, &
D^*\circ D
&=-D\circ D^*=\Ll_{a_0},\\[0.3em]
D(a_0)&=0, &
D^*(a_0)&=0.
\end{aligned}
\end{equation}
Then, the vector space
$
\overline{A} := H^* \oplus A \oplus H$ endowed with the bilinear product $\bullet$ given by
\begin{equation}\label{pro-anti}
\begin{aligned}
d \bullet d &= a_0 + \mu e, & d \bullet u &= D(u) + \omega_A(a_0, u) e, \\
u \bullet d &= -D^*(u) + \omega_A(a_0, u) e, & u \bullet v& = u \bullet_A v + \omega_A(D(u), v) e,
\end{aligned}
\end{equation}
for all $u, v \in A$, is a nearly associative algebra.  

Moreover, the bilinear form $\omega : \overline{A} \times \overline{A} \to \K$ defined by
\[
\omega(u,v) = \omega_A(u,v), \qquad \omega(e,d) =- \omega(d,e) = 1,
\]
and equal to zero otherwise, defines an anti-quadratic  on $(\overline{A}, \bullet)$.
\end{Theorem}

The anti-quadratic nearly associative algebra $(\overline{A}, \bullet, \omega)$ constructed in Theorem~\ref{construction-ex} is called the \emph{double extension} of the anti-quadratic nearly associative algebra $(A, \bullet_A, \omega_A)$ by means of $( D, a_0)$.

\begin{proof}
According to Lemma~\ref{Lemma1}, the relation $D(u \bullet_A v) = -u \bullet_A D(v)$ for all $u, v \in A$ ensures that the bilinear map $\phi(u,v) = \omega_A(D(u), v)$ satisfies Condition~\eqref{line-central}. Hence, the product on $\widetilde{A} := A \oplus H^*$ given by
\[
(u + \alpha e) \widetilde{\bullet} (v + \beta e) = u \bullet_A v + \phi(u,v)e
\]
defines a nearly associative algebra, which corresponds to an annihilator extension of $(A, \bullet_A)$ by $H^*$ induced by of $\phi$.  

Now, since $(\xi, D, a_0, \mu)$ satisfies System~\eqref{con-ex2}, it follows from Lemma~\ref{Lemma2} that $(H, \widetilde{A} := A \oplus H^*, \widetilde{D}, \widetilde{\xi}, \widetilde{a}_0, 0)$ constitutes a context of an almost semi-direct product of nearly associative algebras. Consequently, the generalized semi-direct product $(\overline{A}, \bullet)$ is indeed a nearly associative algebra. 

Furthermore, a straightforward verification shows that the form $\omega$ is anti-invariant on $(\overline{A}, \bullet)$. This implies that $(\overline{A}, \bullet, \omega)$ is an anti-quadratic nearly associative algebra, which completes the proof.
\end{proof}

\sssbegin{Proposition}\label{reduit}
Let $(A, \bullet_A, \omega_A)$ be an anti-quadratic nearly associative algebra. Let $I := \K e \subset N(A)$ be a non-zero vector subspace. We denote by $B := I^\perp / I$ the quotient space and by $\pi : I^\perp \to B$ the canonical projection. Then:
\begin{enumerate}
    \item[$(1)$] $I$ is a totally isotropic two-sided ideal of $(A, \bullet_A)$, and $I^\perp$ is also a two-sided ideal of $(A, \bullet_A)$.
    
    \item[$(2)$] The quotient $B := I^\perp / I$ inherits a canonical structure of an anti-quadratic nearly associative algebra, given by
    \begin{equation}
    \pi(u) \bullet_B \pi(v) = \pi(u \bullet_A v), \qquad \omega_B\big(\pi(u), \pi(v)\big) = \omega_A(u, v),
\label{Bquocient}\end{equation}
    for all $u, v \in I^\perp$.
\end{enumerate}
\end{Proposition}

\begin{proof}
Let us show Part (1). Since $I = \K e \subset N(A)$, we have $I \bullet_A A = \{0\} = A \bullet_A I$, which shows that $I$ is a two-sided ideal. Furthermore, since $I \subset I^\perp$, it is totally isotropic. Now, let $v \in I^\perp$, $w \in A$, and $u \in I$. By the anti-invariance of $\omega_A$, we have
\[
\omega_A(v \bullet_A w, u) = -\omega_A(w, v \bullet_A u) = 0, \qquad \omega_A(w \bullet_A v, u) = -\omega_A(v, w \bullet_A u) = 0,
\]
due to the fact that $v \bullet_A u = w \bullet_A u = 0$. Hence, $v \bullet_A w \in I^\perp$ and $w \bullet_A v \in I^\perp$, which means that $I^\perp$ is a two-sided ideal of $(A, \bullet_A)$.

Let us prove Part (2). Since $I$ and $I^\perp$ are two-sided ideal, the product $\bullet_B$ on the quotient space is well-defined. Thus, $(B, \bullet_B)$ is nearly associative algebra. Moreover, for any $u, v \in I^\perp$ and $x, y \in I$, we have:
\[
\omega_B(u + x, v + y)   = \omega_A(u, v),
\]
since $I \subset I^\perp$. This ensures that the induced bilinear form $\omega_B$ is well-defined and non-degenerate on $B$. The fact that $ \omega_B$ is a  anti-invariant on $(B, \bullet_B)$ follows immediately from the anti-invariance property of $\omega_A$ on $A$. It follows that $(B, \bullet_B, \omega_B)$ is an  anti-quadratic nearly associative algebra.
\end{proof}

We now turn to the converse of Theorem~\ref{construction-ex}, which provides a structural characterization of the algebras obtained by this construction.

\sssbegin{Theorem}\label{Th-ex1}
Let $(A, \bullet_A, \omega_A)$ be an anti-quadratic nearly associative algebra.  There exists a one-dimensional  two-sided ideal $I \subset N(A)$. Moreover, $(A, \bullet_A, \omega_A)$ is a double extension of the anti-quadratic nearly associative algebra 
$
\big(B := I^\perp / I, \, \bullet_B, \, \omega_B\big)
$
by means of $( D, a_0)$.
\end{Theorem}

\begin{proof}
According to Proposition~\ref{nil}, the algebra $(A, \bullet)$ is nilpotent and hence $N(A) \neq \{0\}$.  Let $I := \K e$ be a non-zero vector subspace of $N(A)$. Then, by Proposition~\ref{reduit}, $I$ is a totally isotropic two-sided ideal of $(A,\bullet)$, and its orthogonal  $I^\perp$ is  a two-sided ideal of $(A,\bullet)$. Since $\omega$ is non-degenerate, there exists an element $d \in A \setminus \{0\}$ such that $\omega(e,d) = 1$.

Let $V := \K d$ and $H := (I \oplus V)^\perp$. We have the decomposition
\(
A = I \oplus H \oplus V,\) and \( I^\perp = I \oplus H.
\) For all $u,v \in H$, we can write
\[
u \bullet v = u \bullet_H v + \varphi(u,v)e,
\]  
where $\varphi : H \times H \to \K$ is a bilinear form and $\bullet_H : H \times H \to H$ is a bilinear map.  
It is straightforward to verify that $(H, \bullet_H)$ is a nearly associative algebra and that the restriction
\(
\omega_H :=\omega|_{H \times H}
\)
defines an anti-quadratic structure on $(H,\bullet_H)$.  Therefore, $\varphi$ satisfies Condition~\eqref{line-central}.

According to Proposition~\ref{reduit}, the quotient algebra $(B := I^\perp / I, \, \bullet_B, \, \omega_B)$,
given by \eqref{Bquocient}, is an anti-quadratic nearly associative algebra. The canonical projection $\pi: I^\perp \to B$ allows us to identify $(H,\bullet_H,\omega_H)$ with $(B,\bullet_B,\omega_B)$. Thus, we can write $I^\perp \simeq I \oplus B$. Since $I^\perp$ is a two-sided ideal of $(A,\bullet)$, we obtain
\begin{equation}
d \bullet d = a_0 + \alpha e + \lambda d, \quad 
d \bullet u = D(u) + f(u)e, \quad
u \bullet d = \xi(u) + g(u)e,
\label{pro-anti0}
\end{equation}
for all $u \in B$, where $D, \xi \in \mathrm{End}(B)$, $a_0 \in B$,  $f,g \in B^*$, and   $\alpha, \lambda \in \K$.

Since $\omega$ is anti-quadratic on $(A, \bullet)$ and identifying $H$ with $B$ we obtain the following relations: For all $u, v \in B$, we have
$$f(u) = -\omega(d, d\bullet u) = \omega(d\bullet d, u) = \omega_B(a_0, u),$$
and 
$$g(u) = \omega(d, u\bullet d) = \omega(d\bullet d, u) = \omega_B(a_0, u).$$
Hence, $f(u) = g(u) = \omega_B(a_0, u)$.

Similarly, we have
$$\lambda = \omega(d\bullet d, e) = -\omega(d, d\bullet e) = 0.$$

Moreover, for all $u,v\in B$, we have
$$\varphi(u,v) = -\omega(d, u \bullet v) = \omega(d\bullet u,v) = \omega_B(D(u), v),$$
and 
$$\omega_B(\delta(u), v) = \omega(u\bullet d,v) = -\omega(u, d \bullet v) = -\omega_B(u, D(v)).$$
Consequently, $\varphi(u,v) = \omega_B(D(u), v)$ and $\delta = -D^*$.

Therefore, the product $\bullet$ takes the folloiwng simplified form:
\[
\begin{aligned}
d \bullet d &= a_0 + \mu e, & d \bullet u &= D(u) + \omega_B(a_0, u) e, \\
u \bullet d &= -D^*(u) + \omega_B(a_0, u) e, & u \bullet v& = u \bullet_B v + \omega_B(D(u), v) e,
\end{aligned}
\]
Now, it is straightforward to show that $(A, \bullet)$ is a nearly associative algebra if and only if the triple $(D, a_0)$ satisfies system \eqref{con-ex2}. Hence, $(A, \bullet, \omega)$ is the double extension of the anti-quadratic nearly associative algebra $(B, \bullet_B, \omega_B)$ by means of $(D, a_0)$.
\end{proof}

\sssbegin{Corollary}\label{succession-double}
Any anti-quadratic nearly associative algebra can be obtained by a finite sequence of successive double extensions starting from the trivial algebra $\{0\}$.
\end{Corollary}

\begin{proof}
According to Theorem~\ref{Th-ex1}, $(A, \bullet, \omega)$ is a double extension of a lower-dimensional anti-quadratic nearly associative algebra $(B_1, \bullet_{B_1}, \omega_{B_1})$ by means of $(D_1, (a_0)_1)$. Applying this theorem inductively, $(B_1, \bullet_{B_1}, \omega_{B_1})$ is itself a double extension of a lower-dimensional anti-quadratic nearly associative algebra $(B_2, \bullet_{B_2}, \omega_{B_2})$ by means of $(D_2, (a_0)_2)$. Since $A$ is finite-dimensional, this process of successive double extensions terminates after a finite number of steps at the trivial algebra $\{0\}$.
\end{proof}

\section{Low-dimensional classification of anti-quadratic nearly associative algebra}\label{S6}
In this section, we apply Theorems~\ref{Th-ex1}  to classify all anti-quadratic  nearly associative algebras of dimension up to $4$.

\ssbegin{Definition}[Isometric isomorphism]
Let $(A_1, \bullet_1, \omega_1)$ and $(A_2, \bullet_2, \omega_2)$ be two anti-quadratic algebras. An \emph{isometric isomorphism} between $A_1$ and $A_2$ is an isomorphism of algebras $\phi: A_1 \to A_2$ that preserves the bilinear forms, meaning that:
\[
\omega_2\big(\phi(u), \phi(v)\big) = \omega_1(u, v) \quad \text{for all } u, v \in A_1.
\]
In this case, $(A_1, \bullet_1, \omega_1)$ and $(A_2, \bullet_2, \omega_2)$ are said to be \emph{isometrically isomorphic}.
\end{Definition}

We begin with the classification of anti-quadratic nearly associative algebras of dimension 2.
\ssbegin{Proposition}\label{mil3}
Let $(A, \bullet, \omega)$ be a two-dimensional anti-quadratic nearly associative algebra. Then $(A, \bullet, \omega)$ is is isometrically isomorphic  to one of the following mutually non-isomorphic  anti-quadratic algebras:
\begin{enumerate}
    \item[$(1)$] $(A_1, \omega)$ with the trivial product and $\omega = e_{1}^{*} \wedge e_{2}^{*}$.
    \item[$(2)$] $(A_2,\omega)$: \quad   $e_1\bullet e_1=e_2, \quad   \omega=-e_{1}^{*}\wedge e_{2}^{*}$.
\end{enumerate}
\end{Proposition}

\begin{proof}
Assume that $(A,\bullet,\omega)$ is a non-trivial two-dimensional anti-quadratic nearly associative algebra. By Theorem \ref{Th-ex1}, $(A,\bullet,\omega)$ is a double extension of the zero algebra $\{0\}$ by means of $(0,0,\mu)$, where $\mu\neq0$. It follows from \eqref{pro-anti} that there exists a basis $\{d,e\}$ of $A$ such that
$$
d\bullet d=\mu e, \esp 
\omega=e^{*}\wedge d^{*}.
$$
Setting
$
e_1=\frac{1}{\sqrt[3]{\mu}}d,$ and  $e_2=\sqrt[3]{\mu}e,
$
we obtain
$$
e_1\bullet e_1=e_2, \esp \omega=-e_1^{*}\wedge e_2^{*} .$$
Therefore, $(A,\bullet,\omega)$ is is isometrically isomorphic to $(A_2,\omega)$.
\end{proof}

We now proceed to the classification of anti-quadratic nearly associative algebras of dimension 4. By Theorem~\ref{Th-ex1}, every such algebra can be obtained from a lower-dimensional anti-quadratic nearly associative algebra through a sequence of double extensions.

Let $(A,\bullet,\omega)$ be a $4$-dimensional anti-quadratic nearly associative algebra obtained as a double extension of a $2$-dimensional anti-quadratic nearly associative algebra $(B,\bullet_B,\omega_B)$ by means of a datum $(D, a_0)$. According to Proposition~\ref{mil3}, the algebra $(B,\bullet_B,\omega_B)$ is  isometrically isomorphic  to either the trivial anti-quadratic algebra or to $(A_2,\omega)$.

We first consider the case where $B$ is trivial. In this situation, the conditions of Theorem~\ref{construction-ex}  reduce to
\begin{equation}\label{abelian}
D^2=D^*\circ D=D\circ D^*=0,\qquad
D(a_0)=D^*(a_0)=0.
\end{equation}

\ssbegin{Proposition}\label{Trivial}
Let $(B, \omega_B)$ be a $2$-dimensional  symplectic vector space. Then $(D, a_0)$ satisfies the equations \eqref{abelian} if and only if there exists a symplectic basis $\mathbb{B} = \{e_1, e_2\}$ of $B$ such that the endomorphism $D$ and the element $a_0$ are expressed with respect to $\mathbb{B} $ by one of the following mutually exclusive cases:
\begin{enumerate}
\item[\rm (a)] $D = 0$, and $a_0 = \alpha e_1 + \beta e_2$, where $\alpha, \beta \in \mathbb{K}$. 
\item[\rm (b)] $D = \begin{pmatrix} 0 & x \\ 0 & 0 \end{pmatrix}$, and  $a_0 = \alpha e_1$, where $x \in \mathbb{K}^*$ and $\alpha \in \mathbb{K}$. 
\end{enumerate}
\end{Proposition}

\begin{proof}
Assume that $(D, a_0)$ satisfies  \eqref{abelian}. The condition 
$$
\omega_{B}(D(u), D(v)) = \omega_{B}(u, D^*(D(v))) = 0 \text{  for all $u,v \in \B$,}
$$ implies that $\mathrm{Im}(D)$ is a totally isotropic subspace of $B$. Since $\dim(B) = 2$, we deduce that $\dim(\mathrm{Im}(D)) \leq 1.$ If $\mathrm{Im}(D) = \{0\}$, then $D=0$, which yields case $(a)$.

Now, suppose that $\mathrm{Im}(D) \neq \{0\}$, implying $\dim(\mathrm{Im}(D)) = 1$. We can then choose a symplectic basis $\{e_1, e_2\}$ of $(B, \omega_{B})$ such that the matrix matrix representation of $D$ has the form:
\[
D = \begin{pmatrix} y & x \\ 0 & 0 \end{pmatrix}, \quad \text{with } (y, x) \neq (0,0).
\]

From  $D^2 = 0$, a direct calculation  shows that $y=0$. Consequently, condition $D \circ D^* = D^* \circ D = 0$ is automatically satisfied.
Finally, conditions $D(a_0) = D^*(a_0) = 0$ imply that $a_0 = \alpha e_1$ for some scalar $\alpha \in \mathbb{K}$.
\end{proof}

\ssbegin{Proposition}\label{ntri}
Let $(A_2, \omega)$ be the unique non-trivial anti-quadratic nearly associative  algebra of dimension $2$. Then, $(D, a_0)$ satisfies system \eqref{con-ex2} if and only if there exists a  basis $\{e_1, e_2\}$ of $A_2$ such that $\omega=-e_{1}^{*} \wedge e_{2}^{*}$,  the endomorphism $D$ and the element $a_0$ are expressed by one of the following cases:
\begin{enumerate}
    \item[\rm (a)] $ D = 0$, and $a_0 = \beta e_2$, with $\beta \in \mathbb{K}$.
    
    \item[\rm (b)] $D = \begin{pmatrix} 0 & 0 \\ x & 0 \end{pmatrix}$, and  $a_0 = \beta e_2$, with $x \in \mathbb{K}\setminus \{0\}$ and $\beta\in \mathbb{K}.$
\end{enumerate}
\end{Proposition}

\begin{proof}
By Proposition~\ref{mil3}, there exists a symplectic basis $\{e_1,e_2\}$ of $(A_2,\omega)$ such that
$$
e_1\bullet e_1=e_2,
\qquad
\omega(e_1,e_2)=-1.
$$
If $D=0$, then the relation $\Ll_{a_0}=0$ implies $a_0=\beta e_2$ for some $\beta\in\K$. This yields case $(a)$.

Assume now that $D\neq0$. Let
$
D(e_1)=a_{11}e_1+a_{21}e_2,$ and $D(e_2)=a_{12}e_1+a_{22}e_2.$ Using the identity $
D(u\bullet v)=-u\bullet D(v),$
we obtain
$
D(e_2)=D(e_1\bullet e_1)
=-e_1\bullet D(e_1)
=-a_{11}e_2.
$
Consequently,
$
a_{12}=0,$ and $a_{22}=-a_{11}$. Since $D^2=0$, we get
$
0=D^2(e_2)=D(-a_{11}e_2)=a_{11}^{2}e_2,$
and therefore $a_{11}=0$. Hence
$D(e_1)=a_{21}e_2,$ and 
$D(e_2)=0$. Since $D\neq 0$, then $a_{21}\neq0$. The condition $D(a_0)=0$ then implies that $a_0=\beta e_2$ for some $\beta\in\K$. Finally, a direct verification shows that all the remaining relations in \eqref{con-ex2} are automatically satisfied. Therefore, we obtain case $(b)$.
\end{proof}

We have now all the needed ingredients to give the complete classification list of four-dimensional anti-quadratic nearly associative  algebras.

\ssbegin{Theorem}\label{classif4}
Every  $4$-dimensional anti-quadratic nearly associative algebra  $(A, \bullet, \omega)$ is  isometrically isomorphic  to one of the anti-quadratic nearly associative algebras:
\begin{enumerate}
    \item[\rm (1)] $(A_{4,1}, \omega)$ with the trivial product and $\omega = e_1^* \wedge e_2^* + e_3^* \wedge e_4^*$.
    \item[\rm (2)] $A_{4,2} = \big(\mathrm{span}(d, e_1, e_2, e), \bullet, \omega\big)$ where:
    \[
    \begin{cases} 
    d\bullet d=\al e_1+\beta e_2+\mu e,\; d\bullet e_1=e_1\bullet d=-\beta e,\; d\bullet e_2=e_2\bullet d=\al e, \\
    \omega = e^* \wedge d^* + e_1^* \wedge e_2^*, \quad \text{with } \alpha, \beta, \mu \in \mathbb{K}.
    \end{cases}
    \]
 \item[\rm (3)] $A_{4,3} = \big(\mathrm{span}(d, e_1, e_2, e), \bullet, \omega\big)$ where:
    \[
    \begin{cases} 
    d\bullet d=\al e_1+\mu e,\;\; d\bullet e_2=e_2\bullet d=xe_1+\al e,\; e_2\bullet e_2=x e, \\
    \omega = e^* \wedge d^* + e_1^* \wedge e_2^*, \quad \text{with } x\in \mathbb{K}^* , \al, \mu \in \mathbb{K}.
    \end{cases}
    \]

\item[\rm (4)] $A_{4,4} = \big(\mathrm{span}(d, e_1, e_2, e), \bullet, \omega\big)$ where: 
    \[
    \begin{cases} 
    d\bullet d=\beta e_2+\mu e,\;\; d\bullet e_1=e_1\bullet d=\beta  e, \; e_1\bullet e_1=e_2,\\
    \omega = e^* \wedge d^* - e_1^* \wedge e_2^*, \quad \text{with } \beta, \mu \in \mathbb{K}.
    \end{cases}
    \]

\item[\rm (5)] $A_{4,5} = \big(\mathrm{span}(d, e_1, e_2, e), \bullet, \omega\big)$ where: 
\[
    \begin{cases} 
    d\bullet d=\beta e_2+\mu e,\;\; d\bullet e_1=e_1\bullet d=xe_2+\beta  e, \; e_1\bullet e_1=e_2+x e,\\
    \omega = e^* \wedge d^* - e_1^* \wedge e_2^*, \quad \text{with }  x\in \mathbb{K}\setminus \{0\}, \beta, \mu \in \mathbb{K}.
    \end{cases}
    \]
    \end{enumerate}
    \end{Theorem}
   
\begin{proof}
According to Theorem \ref{Th-ex1}, every $4$-dimensional anti-quadratic nearly associative algebra $(A, \bullet, \omega)$ can be realized as a double extension of a $2$-dimensional anti-quadratic nearly associative algebra $(B, \bullet_B, \omega_B)$ by means of  $(D, a_0)$. 

By virtue of Proposition \ref{Trivial} and Proposition \ref{ntri}, the solutions for the structural  $(D, a_0)$ are explicitly determined by distinguishing two cases depending on whether the algebra $B$ is trivial or non-trivial. Finally, applying the general product formula given in \eqref{pro-anti} to each of these cases yields the desired classification of the anti-quadratic algebras listed above.
\end{proof}
\section*{Acknowledgements}
The authors thank Sa\"id Benayadi for pointing out references \cite{Lyla1, Layla3}.
\section*{Data Availability}
Data sharing is not applicable to this article, as no datasets were
generated or analyzed during the current study.
\section*{Author Contributions}
All authors contributed equally to this work and approved the final
manuscript.
\section*{Competing Interests}
The authors declare that they have no competing interests.
\bibliographystyle{elsarticle-num}

\end{document}